\documentclass[10pt, a4paper, fleqn]{article}
\usepackage{titling}
\title{Projective and Telescopic Projective Integration for Non-Linear Kinetic Mixtures}

\newcommand{\authorPDF}{Bailo, Rey.}
\newcommand{\subjectPDF}{65M08; 65M12; 76P05; 82C40. }
\newcommand{\keywordsPDF}{Multispecies gas; kinetic mixture; Boltzmann equation; BGK model; projective integration; Sod Tube; Kelvin-Helmholtz; Richtmyer-Meshkov.}
 \usepackage{authblk}

\author[1,2]{Rafael Bailo}
\author[1]{Thomas Rey}

\affil[1]{
Univ. Lille, CNRS, Inria, UMR 8524 - Laboratoire Paul Painlev\'{e}
}
\affil[ ]{
	F-59000 Lille, France
}
\affil[ ]{\textit{
		rafael.bailo@univ-lille.fr
	}}
\affil[ ]{\textit{
		thomas.rey@univ-lille.fr
	}}

\affil[ ]{}

\affil[2]{
	Mathematical Institute, University of Oxford
}
\affil[ ]{
	OX2 6GG Oxford, United Kingdom
}
\affil[ ]{\textit{bailo@maths.ox.ac.uk}} 
\usepackage[left=1.1in, right=1.1in, top=1.1in, bottom=1.1in]{geometry}

\makeatletter
\let\newtitle\@title
\let\newauthor\@author
\let\newdate\@date
\makeatother
\usepackage{fancyhdr}
\usepackage{setspace}
\usepackage[
	hypertexnames=false,
	pdftex,
	pdftitle={\thetitle. \authorPDF},
	pdfauthor={\authorPDF},
	pdfsubject={\subjectPDF},
	pdfkeywords={\keywordsPDF},
]{hyperref}
\hypersetup{
	colorlinks=true,
	linkcolor=blue,
	citecolor=blue,
	urlcolor=blue,
	linktocpage
}

\usepackage[utf8]{inputenc}
\usepackage[T1]{fontenc}
\usepackage{amsmath}
\usepackage{amsthm}
\usepackage{amsfonts}
\usepackage{amssymb}
\usepackage{graphicx}
\usepackage{mathtools}
\usepackage{ifthen}
\usepackage{dsfont}

\usepackage[UKenglish]{babel}

\usepackage[dvipsnames]{xcolor}

\definecolor{ppGreen}{HTML}{008000}
\definecolor{ppBlue}{HTML}{0000FF}
\definecolor{ppRed}{HTML}{FF0000}
\definecolor{ppPurple}{HTML}{800080}
\definecolor{lightblue}{rgb}{0.145,0.6666,1}
\definecolor{grey52}{RGB}{52,52,52}
\definecolor{color1}{RGB}{0,62,116}
\definecolor{color2}{RGB}{152,152,152}
\definecolor{color3}{RGB}{52,52,52}
\definecolor{color4}{RGB}{100,100,100}

\definecolor{imperialnavy}{RGB}{0,33,71}
\definecolor{imperialblue}{RGB}{0,62,116}
\definecolor{imperialgrey}{RGB}{235,238,238}
\definecolor{imperialcoolgrey}{RGB}{157,157,157}

\definecolor{lille}{RGB}{178, 35, 114} 

\newcounter{review}

\newcommand{\ntcreview}[3]{\refstepcounter{review}

	{\color{#2}{\textbf{[#1]}: #3}}}

\newcommand{\creview}[3]{\ntcreview{#1}{#2}{#3}
	\addcontentsline{tor}{subsection}{\thereview~\textbf{[#1]}:~#3
	}}

\newcommand{\review}[2]{\creview{#1}{blue}{#2}}

\makeatletter

\newcommand\listreviewname{List of Reviews}
\newcommand\listofreviews{\section*{\listreviewname}\@starttoc{tor}}
\makeatother

\usepackage[all]{nowidow}

\newcommand{\subjectclassification}[1]{

	{\small\textbf{\textit{AMS Subject Classification --- }} #1}

}
\newcommand{\keywords}[1]{

	{\small\textbf{\textit{Keywords --- }} #1}

}

\renewcommand\lll\MoveEqLeft

\usepackage{tikz}
\usetikzlibrary{patterns}
\usetikzlibrary{decorations.pathreplacing}
\usetikzlibrary{calc}
\usetikzlibrary{fadings}
\usetikzlibrary{shapes}
\usetikzlibrary{plotmarks}
\usetikzlibrary{arrows, decorations.markings}
\tikzset{thicker line small arrows m/.style args={#1in#2}{
			draw=#2,
			solid,
			line width=#1,
			shorten >=1mm,
			decoration={
					markings,
					mark=at position 1.0 with {\arrow[fill=#2,thin]{triangle 90}}
				},
			postaction={decorate}
		}}

\usepackage{pgfplots}
\pgfplotsset{compat=1.16}

\usepackage{float}
\usepackage{subcaption}
\usepackage[section]{placeins}
\usepackage{rotating}

\usepackage{array}
\newcolumntype{L}[1]{>{\raggedright\let\newline\\\arraybackslash\hspace{0pt}}m{#1}}
\newcolumntype{C}[1]{>{\centering\let\newline\\\arraybackslash\hspace{0pt}}m{#1}}
\newcolumntype{R}[1]{>{\raggedleft\let\newline\\\arraybackslash\hspace{0pt}}m{#1}}
\usepackage{booktabs}
\usepackage{tabularx}
\usepackage{multirow}

\usepackage{titling}

\usepackage[algoruled]{algorithm2e}

\usepackage[capitalise]{cleveref}

 \usepackage{accents}

\newcommand\term\emph

\numberwithin{equation}{section}

\usepackage{autobreak}

\usepackage{enumerate}

\newcommand{\singleappendixtitle}[1]{
	\appendix

	\section*{Appendix: #1}

	\addcontentsline{toc}{section}{Appendix}

	\setcounter{equation}{0}
	\renewcommand\theequation{A.\arabic{equation}}

	\setcounter{theorem}{0}
	\renewcommand\thetheorem{A.\arabic{theorem}}
}

\makeatletter
\def\@maketitle{\newpage
	\begin{center}\let \footnote \thanks
		{\LARGE\bfseries \@title \par}\vskip 2.5em{\large
				\lineskip .5em\begin{tabular}[t]{c}\@author
				\end{tabular}\par}\vskip 1em{\large \@date}\end{center}\par
	\vskip 1.5em}
\makeatother

\newenvironment{system}{
	\begingroup
	\left\lbrace \begin{array}{{>{\displaystyle\vphantom{\frac{1}{1}}}l}}
		}{
	\end{array} \right.
	\endgroup
}



\usepackage{amsthm}
\theoremstyle{plain}

\theoremstyle{remark}

\usepackage{esint}

\def\XXint#1#2#3{{\setbox0=\hbox{$#1{#2#3}{\int}$ }
			\vcenter{\hbox{$#2#3$ }}\kern-.6\wd0}}

\renewcommand{\th}{\textsuperscript{th} }

\DeclarePairedDelimiter{\prt}{(}{)}
\DeclarePairedDelimiter{\brk}{[}{]}

\DeclarePairedDelimiter{\abs}{|}{|}
\DeclarePairedDelimiter{\norm}{\|}{\|}
\DeclarePairedDelimiter{\set}{\{}{\}}

\DeclarePairedDelimiter{\inn}{\langle}{\rangle}

\usepackage{suffix}

\newcommand{\inner}[2]{\inn{#1,#2}}
\WithSuffix\newcommand\inner*[2]{\inn*{#1,#2}}

\DeclarePairedDelimiter{\positive}{(}{)^{+}}
\DeclarePairedDelimiter{\negative}{(}{)^{-}}

\newcommand\pos\positive
\renewcommand\neg\negative

\WithSuffix\newcommand\pos*{\positive*}
\WithSuffix\newcommand\neg*{\negative*}

\newcommand{\R}{{\mathbb{R}}}

\renewcommand{\L}[1]{{L^{#1}}}
\newcommand{\Lone}{\L{1}}
\newcommand{\Ltwo}{\L{2}}

\newcommand{\curlyX}{\mathcal{X}}
\newcommand{\curlyV}{\mathcal{V}}

\newcommand{\pnorm}[2]{\norm{#2}_{\L{#1}}}
\WithSuffix\newcommand\pnorm*[2]{\norm*{#2}_{\L{#1}}}

\newcommand{\psnorm}[3]{\norm{#3}_{\L{#1}(#2)}}
\WithSuffix\newcommand\psnorm*[3]{\norm*{#3}_{\L{#1}(#2)}}

\newcommand{\pnormp}[2]{\pnorm{#1}{#2}^{#1}}
\WithSuffix\newcommand\pnormp*[2]{\pnorm*{#1}{#2}^{#1}}

\newcommand{\psnormp}[3]{\psnorm{#1}{#2}{#3}^{#1}}
\WithSuffix\newcommand\psnormp*[3]{\psnorm*{#1}{#2}{#3}^{#1}}

\newcommand\svec\vec

\renewcommand{\vec}{\mathbf}
\renewcommand{\svec}{\boldsymbol}

\newcommand{\bx}{\vec{x}}

\newcommand{\bv}{\vec{v}}

\renewcommand{\d}{\mathrm{d}}
\newcommand{\dd}{\mathop{}\!\d}

\newcommand{\pder}[2]{\frac{\partial #1}{\partial #2}}

\newcommand{\dbx}{\dd \bx}

\newcommand{\grad}{\nabla}

\newcommand{\pt}{\partial_t}

\newcommand{\Dt}{\Delta t}
\newcommand{\Dx}{\Delta x}
\newcommand{\Dy}{\Delta y}

\newcommand{\nhalf}{1/2}

\renewcommand{\i}{_{i}}

\newcommand{\ih}{_{i+\nhalf}}
\newcommand{\imh}{_{i-\nhalf}}

\renewcommand{\j}{_{j}}

\newcommand{\n}{^{n}}
\newcommand{\np}{^{n+1}}

\newcommand{\ppr}{(r)}

\newcommand{\Wr}{^{W,\,\ppr}}

\newlength{\dhatheight}

\ifthenelse{\isundefined{\Wr}}{
	\newcommand{\Wr}{^{W,\,\ppr}}
}{
	\renewcommand{\Wr}{^{W,\,\ppr}}
}

\newcommand{\curlyN}{\mathcal{N}}

 \usepackage{todonotes}

\newif\ifskiptable

\newcommand{\placedfigure}[1]{
	\begin{figure}\centering
		#1
	\end{figure}
}

\newcommand{\placedfigureHT}[1]{
	\begin{figure}[!ht]
		\centering
		#1
	\end{figure}
}

\pgfplotsset{colormap={hsv}{
			hsb(0.00cm)=(0.00,0,0.95);
			hsb(0.05cm)=(0.05,1,1);
			hsb(0.10cm)=(0.10,1,1);
			hsb(0.15cm)=(0.15,1,1);
			hsb(0.20cm)=(0.20,1,1);
			hsb(0.25cm)=(0.25,1,1);
			hsb(0.30cm)=(0.30,1,1);
			hsb(0.35cm)=(0.35,1,1);
			hsb(0.40cm)=(0.40,1,1);
			hsb(0.45cm)=(0.45,1,1);
			hsb(0.50cm)=(0.50,1,1);
			hsb(0.55cm)=(0.55,1,1);
			hsb(0.60cm)=(0.60,1,1);
			hsb(0.65cm)=(0.65,1,1);
			hsb(0.70cm)=(0.70,1,1);
			hsb(0.75cm)=(0.75,1,1);
			hsb(0.80cm)=(0.80,1,1);
			hsb(0.85cm)=(0.85,1,1);
			hsb(0.90cm)=(0.90,1,1);
			hsb(0.95cm)=(0.95,1,1);
			hsb(1.00cm)=(1.00,1,1);
		}
}

\pgfplotsset{colormap={hsvSoft}{
			hsb(0.00cm)=(0.00,0,0.95);
			hsb(0.05cm)=(0.05,1,1);
			hsb(0.10cm)=(0.10,1,1);
			hsb(0.15cm)=(0.15,1,1);
			hsb(0.20cm)=(0.20,1,1);
			hsb(0.25cm)=(0.25,1,1);
			hsb(0.30cm)=(0.30,1,1);
			hsb(0.35cm)=(0.35,1,1);
			hsb(0.40cm)=(0.40,1,1);
			hsb(0.45cm)=(0.45,1,1);
			hsb(0.50cm)=(0.50,1,1);
			hsb(0.55cm)=(0.55,1,1);
			hsb(0.60cm)=(0.60,1,1);
			hsb(0.65cm)=(0.65,1,1);
			hsb(0.70cm)=(0.70,1,1);
			hsb(0.75cm)=(0.75,1,1);
			hsb(0.80cm)=(0.80,1,1);
			hsb(0.85cm)=(0.85,1,1);
			hsb(0.90cm)=(0.90,1,1);
			hsb(0.95cm)=(0.95,1,1);
			hsb(1.00cm)=(0.00,0,0.95);
		}
}

\pgfplotsset{colormap={viridisSoft}{
			rgb255=(242, 242, 242);
			rgb255=(242, 242, 242);
			rgb=(0.28026,0.1657,0.4765);
			rgb=(0.26366,0.23763,0.51877);
			rgb=(0.23744,0.3052,0.54192);
			rgb=(0.20862,0.36775,0.55267);
			rgb=(0.18225,0.42618,0.55711);
			rgb=(0.1592,0.48224,0.55807);
			rgb=(0.13777,0.53749,0.5549);
			rgb=(0.12115,0.59274,0.54465);
			rgb=(0.12808,0.64775,0.5235);
			rgb=(0.18065,0.7014,0.48819);
			rgb=(0.27415,0.75198,0.4366);
			rgb=(0.39517,0.79747,0.36775);
			rgb=(0.53561,0.83578,0.2819);
			rgb=(0.68895,0.86545,0.18272);
			rgb=(0.84557,0.88733,0.0997);
			rgb=(0.99324,0.90616,0.14394)
		}
} %

 \newcommand{\KHTwoDPressureShortFig}[6]{

	\begin{tikzpicture}
		\renewcommand*\showkeyslabelformat[1]{}
		\makeatletter
		\def\SK@@ref#1>#2\SK@{}
		\makeatletter

		\begin{groupplot}[
				group style={
						group name=topLeft,
						group size=1 by 1,
						vertical sep=\plotSmallSeparation,
					},
				colormap name = viridis,
				view = {0}{90},
				xlabel={},
				xticklabel=\empty,
				ylabel={},
				yticklabel=\empty,
				width = \mylinewidth,
				height = \mylinewidth,
				colorbar,
				colorbar style={
						at={(-0.15,1)}, anchor=above north west, width=\plotBarWidth,
yticklabel pos=left,
						ylabel={Density $\rho_1$},
					},
				point meta min = 0.0,
				point meta max = 1.05,
				xmin=-0.5,
				xmax=0.5,
				ymin=-0.5,
				ymax=0.5,
			]

			\nextgroupplot[
				title={$t=1.3\times 10^{-2}$},
				title style={yshift=\plotGroupTitleRaise},
				legend to name={commonLegend},
				legend style={
						legend columns=2,
					},
			]

		\end{groupplot}

		\begin{groupplot}[
				group style={
						group name=topRight,
						group size=1 by 1,
						vertical sep=\plotSmallSeparation,
					},
				colormap name = viridis,
				view = {0}{90},
				xlabel={},
				xticklabel=\empty,
				ylabel={$y$},
				yticklabel pos=right,
				width = \mylinewidth,
				height = \mylinewidth,
point meta min = 0.0,
				point meta max = 1.05,
				xmin=-0.5,
				xmax=0.5,
				ymin=-0.5,
				ymax=0.5,
			]

			\nextgroupplot[
				title={$t=2.1\times 10^{-2}$},
				title style={yshift=\plotGroupTitleRaise},
				anchor=north west,
				at={($(topLeft c1r1.north east) + (\plotBigSeparation,0)$)},
			]

		\end{groupplot}

		\begin{groupplot}[
				group style={
						group name=midLeft,
						group size=1 by 1,
						vertical sep=\plotSmallSeparation,
					},
				colormap name = viridis,
				view = {0}{90},
				xlabel={},
				xticklabel=\empty,
				ylabel={},
				yticklabel=\empty,
				width = \mylinewidth,
				height = \mylinewidth,
				colorbar,
				colorbar style={
						at={(-0.15,1)}, anchor=above north west, width=\plotBarWidth,
yticklabel pos=left,
						ylabel={Density $\rho_2$},
					},
				point meta min = 0.0,
				point meta max = 2.15,
				xmin=-0.5,
				xmax=0.5,
				ymin=-0.5,
				ymax=0.5,
			]

			\nextgroupplot[
				anchor=north west,
				at={($(topLeft c1r1.south west) + (0,-\plotBigSeparation)$)},
				title={},
				title style={yshift=\plotGroupTitleRaise},
				legend to name={commonLegend},
				legend style={
						legend columns=2,
					},
			]

		\end{groupplot}

		\begin{groupplot}[
				group style={
						group name=midRight,
						group size=1 by 1,
						vertical sep=\plotSmallSeparation,
					},
				colormap name = viridis,
				view = {0}{90},
				xlabel={},
				xticklabel=\empty,
				ylabel={$y$},
				yticklabel pos=right,
				width = \mylinewidth,
				height = \mylinewidth,
point meta min = 0.0,
				point meta max = 2.15,
				xmin=-0.5,
				xmax=0.5,
				ymin=-0.5,
				ymax=0.5,
			]

			\nextgroupplot[
				title={},
				title style={yshift=\plotGroupTitleRaise},
				anchor=north west,
				at={($(midLeft c1r1.north east) + (\plotBigSeparation,0)$)},
			]

		\end{groupplot}

		\begin{groupplot}[
				group style={
						group name=botLeft,
						group size=1 by 1,
						vertical sep=\plotSmallSeparation,
					},
				colormap name = viridis,
				view = {0}{90},
				xlabel={$x$},
ylabel={},
				yticklabel=\empty,
				width = \mylinewidth,
				height = \mylinewidth,
				colorbar,
				colorbar style={
						at={(-0.15,1)}, anchor=above north west, width=\plotBarWidth,
yticklabel pos=left,
						ylabel={Total Pressure $p$},
					},
				point meta min = 0.5,
				point meta max = 1.15,
				xmin=-0.5,
				xmax=0.5,
				ymin=-0.5,
				ymax=0.5,
			]

			\nextgroupplot[
				anchor=north west,
				at={($(midLeft c1r1.south west) + (0,-\plotBigSeparation)$)},
				title={},
				title style={yshift=\plotGroupTitleRaise},
				legend to name={commonLegend},
				legend style={
						legend columns=2,
					},
			]

		\end{groupplot}

		\begin{groupplot}[
				group style={
						group name=botRight,
						group size=1 by 1,
						vertical sep=\plotSmallSeparation,
					},
				colormap name = viridis,
				view = {0}{90},
				xlabel={$x$},
ylabel={$y$},
				yticklabel pos=right,
				width = \mylinewidth,
				height = \mylinewidth,
				point meta min = 0.5,
				point meta max = 1.15,
				xmin=-0.5,
				xmax=0.5,
				ymin=-0.5,
				ymax=0.5,
			]

			\nextgroupplot[
				title={},
				title style={yshift=\plotGroupTitleRaise},
				anchor=north west,
				at={($(botLeft c1r1.north east) + (\plotBigSeparation,0)$)},
			]

			\addplot3 [
				surf,
				mesh/ordering=y varies,
				shader=flat,
			] table
				{\detokenize{#6}};

		\end{groupplot}

	\end{tikzpicture}

} %
  
\newcommand\eps\varepsilon
\newcommand{\e}{^{\eps}}
\newcommand{\Dv}{\Delta v}
\renewcommand{\O}{\mathcal{O}}
\newcommand{\curlyJ}{\mathcal{J}}
\newcommand{\curlyH}{\mathcal{H}}
\newcommand{\curlyS}{\mathcal{S}}
\newcommand{\curlyP}{\mathcal{P}}

\newcommand{\bbf}{\svec{f}}
\newcommand{\curlyM}{\mathcal{M}}

\newcommand{\dbv}{\dd \bv}
\newcommand{\barv}{\bar{v}}
\newcommand{\barbv}{\bar{\bv}}

\newcommand{\dimx}{{D_\bx}}
\newcommand{\dimv}{{D_\bv}}
\newcommand{\Rx}{{\R^{\dimx}}}
\newcommand{\Rv}{{\R^{\dimv}}}

\newcommand{\Q}{\mathcal{Q}}
\newcommand{\D}{\mathcal{D}}
\newcommand{\M}{\mathcal{M}}

\newcommand{\QBGK}{\Q^\text{BGK}}

\newcommand{\p}{_{p}}
\newcommand{\q}{_{q}}
\newcommand{\pq}{_{p,\,q}}
\newcommand{\qp}{_{q,\,p}}
\newcommand{\pp}{_{p,\,p}}

\newcommand{\divx}{\nabla_\bx\cdot}

\newcommand{\eq}{_{\text{eq}}}

\newcommand{\pij}{_{p,\,i,\,j}}
\newcommand{\pihj}{_{p,\,i+\nhalf,\,j}}
\newcommand{\pimhj}{_{p,\,i-\nhalf,\,j}}

\newcommand{\bj}{{\svec{j}}}
\newcommand{\pbj}{_{p,\,\bj}}
\newcommand{\pqbj}{_{p,\,q,\,\bj}}

\newcommand{\bi}{{\svec{i}}}

\newcommand{\Dbv}{\Delta\bv}

\newcommand{\bC}{\svec{C}}

\newcommand{\pbi}{_{p,\,\bi}}
\newcommand{\pbibj}{_{p,\,\bi,\,\bj}}
\newcommand{\pqbibj}{_{p,\,q,\,\bi,\,\bj}}

\newcommand{\pipj}{_{p,\,i+1,\,j}}

\newcommand{\projnk}{^{n,\,k}}
\newcommand{\projnkp}{^{n,\,k+1}}
\newcommand{\projnzero}{^{n,\,0}}
\newcommand{\projnpzero}{^{n+1,\,0}}

\newcommand{\projnkzero}{^{n,\,k,\,0}}
\newcommand{\projnK}{^{n,\,K}}
\newcommand{\projnKp}{^{n,\,K+1}}
\newcommand{\projnkl}{^{n,\,k,\,l}}
\newcommand{\projnklp}{^{n,\,k,\,l+1}}
\newcommand{\projnkL}{^{n,\,k,\,K_0}}
\newcommand{\projnkLp}{^{n,\,k,\,K_0+1}}

\newcommand{\telnK}{^{n,\,K_1}}
\newcommand{\telnKp}{^{n,\,K_1+1}}

\DeclareMathOperator{\FE}{FE}
\DeclareMathOperator{\PFE}{PFE}

\usepgfplotslibrary{patchplots}

\pgfplotsset{every axis/.append style={
			grid=both,
			grid style={white, line width=.1pt},
			major grid style={white, line width=0.5pt},
			axis background/.style={fill=gray!10},
			axis line style={draw=none},
			tick style={draw=none},
			xlabel = $x$,
line width=1pt,
legend style={
					line width = 1pt,
					draw=none,
					/tikz/every even column/.append style={column sep=0.5cm}
				},
		}}

\newcommand{\plotBigSeparation}{15pt}

\newcommand{\plotSmallSeparation}{6pt}

\newcommand{\plotGroupTitleRaise}{-0.6em} 
\newcommand{\plotBarWidth}{7pt}

\usepgfplotslibrary{groupplots}

\usepackage{calc}

\definecolor{gg0}{HTML}{E24A33}
\definecolor{gg1}{HTML}{348ABD}
\definecolor{gg2}{HTML}{988ED5}
\definecolor{gg3}{HTML}{777777}
\definecolor{gg4}{HTML}{FBC15E}
\definecolor{gg5}{HTML}{8EBA42}
\definecolor{gg6}{HTML}{FFB5B8}

\pgfplotsset{
	/pgfplots/colormap={bright}{rgb255=(0,0,0) rgb255=(78,3,100) rgb255=(2,74,255)
			rgb255=(255,21,181) rgb255=(255,113,26) rgb255=(147,213,114) rgb255=(230,255,0)
			rgb255=(255,255,255)}
}

\newcommand\Mach{\mbox{\text{Ma}}}

\usepackage{siunitx}
\usepackage{chemformula}

\usepackage[nolist]{acronym}  

\renewcommand{\review}[2]{}
\renewcommand{\creview}[3]{}
\renewcommand{\ntcreview}[3]{}

\renewcommand{\tableofcontents}{}
\renewcommand{\listofreviews}{}

\expandafter\def\csname ver@etex.sty\endcsname{3000/12/31}

\usepackage{autonum}
\makeatletter
\patchcmd{\autonum@saveEnvironmentSubcommands}
{(0,0)\begin}
{(0,0)\hfuzz=\maxdimen\begin}
{}{}
\makeatother

\AtBeginDocument{

}

\makeatletter
\autonum@generatePatchedReferenceCSL{ref}
\autonum@generatePatchedReferenceCSL{eqref}
\autonum@generatePatchedReferenceCSL{Cref}
\makeatother

\newcommand*\showkeyslabelformat[1]{} 

\newcommand{\revision}[2]{#2}
\newcommand{\revisionNote}[2]{}

\definecolor{revisionColourOne}{RGB}{180,0,0}
\definecolor{revisionColourTwo}{RGB}{102,0,51}

\newcommand{\revisionOne}[1]{\revision{revisionColourOne}{#1}}
\newcommand{\revisionTwo}[1]{\revision{revisionColourTwo}{#1}}

\newcommand{\revisionNoteOne}[1]{\revisionNote{revisionColourOne}{#1}}
\newcommand{\revisionNoteTwo}[1]{\revisionNote{revisionColourTwo}{#1}}

\begin{document}

\begin{singlespace}\maketitle\end{singlespace}
\begin{abstract}
	We propose fully explicit projective integration and telescopic projective integration schemes for the multispecies Boltzmann and \acf{BGK} equations. The methods employ a sequence of small forward-Euler steps, intercalated with large extrapolation steps. The telescopic approach repeats said extrapolations as the basis for an even larger step. This hierarchy renders the computational complexity of the method essentially independent of the stiffness of the problem, which permits the efficient solution of equations in the hyperbolic scaling with very small Knudsen numbers. We validate the schemes on a range of scenarios, demonstrating its prowess in dealing with extreme mass ratios, fluid instabilities, and other complex phenomena.
\end{abstract}

\subjectclassification{\subjectPDF}
 \revisionNoteOne{Changes corresponding to the comments of Reviewer 1 are shown in red.}
\revisionNoteTwo{Similarly, changes corresponding to Reviewer 2 are shown in purple.}
\tableofcontents
\listofreviews

\begin{acronym}
	\acro{AP}[AP]{asymptotic preserving}
	\acro{BGK}[BGK]{Bhatnagar-Gross-Krook}
	\acro{PI}[PInt]{projective integration}
	\acro{TPI}[TPInt]{telescopic projective integration}
\end{acronym} \section{Introduction}

Mixtures of rarefied gases are found in a wide variety of systems, ranging from the re-entry of an interplanetary probe in the upper atmosphere \cite{D.C2008} to microscale flows in pumps which use no moving parts (viz. the Knudsen compressor, \cite{K.T.O+2004}). Such complex gases cannot be described by classical fluid models, such as the compressible Euler or the Navier-Stokes systems, because of their non-equilibrium behaviour, induced by their rarefaction. Kinetic models, such as the seminal Boltzmann equation, are thus favoured to describe these systems because they are able to reflect the non-equilibrium character of the gases, retaining information about the microscopic many-particle dynamics while avoiding the sheer complexity of the microscopic approach.
Furthermore, actual gases are usually a mixture of different chemical species; for instance, the chemistry of the upper atmosphere is made of up to 20 different species (mostly recombinations of \ch{O2}, \ch{CO2}, \ch{H2}, and \ch{CH4}). The realistic simulations of such systems must involve multispecies kinetic models.

Whereas the mathematical properties of the classical Boltzmann equation for a single species gas are well known by now (crucially, its derivation from Newtonian dynamics was addressed in \cite{G.S.T2013}), many questions remain open for the multispecies case. The most recent theoretical results on the topic come from the series of papers \cite{Briant2016,B.D2016}; they prove the existence, uniqueness, positivity, and exponential trend to equilibrium for the full non-linear multispecies Boltzmann equation in a perturbative, polynomially-weighted, and isotropic $L^1_v L^\infty_x$ setting. The case of the unperturbed setting remains mostly open.

At the numerical level, the most advanced deterministic methods use a Fourier approach to evaluate the full collision operator: a fast spectral algorithm was recently introduced in \cite{W.Z.R+2015} for the multispecies case. Nevertheless, the method was not \ac{AP}, namely, stable in the small relaxation parameter limit; the only recent paper on this subject can be found in \cite{B.B.G2019}, where the \revisionTwo{Maxwell-Stefan} limit for a multispecies gas is investigated numerically via a moment method (which does not compute the full operator). In this work, we shall introduce a new family of numerical integrators for full kinetic multispecies models that are able to deal with a large range of values of the relaxation parameters uniformly on the numerical parameters.

Our work will rely on \ac{PI}, a robust and fully explicit method that allows for the time integration of (two-scale) stiff systems with arbitrary order of accuracy in time. The method was first proposed in \cite{G.K2003} for stiff systems of ordinary differential equations with a clear gap in their eigenvalue spectrum. In such problems, the fast modes, corresponding to the Jacobian eigenvalues with large negative real parts, decay quickly; it is the slow modes, related to eigenvalues of smaller magnitude, that are of practical interest. In this regime, \ac{PI} permits a stable yet explicit integration by combining small and large steps. The integrator performs a few small (inner) steps of an explicit method, using a step size $\delta t$ , until the transients corresponding to the fast modes have decayed; subsequently, the solution is projected (extrapolated) forward in time over a large (outer) time step of size ${\Delta t \gg \delta t}$.

\Ac{PI} was analysed for kinetic equations with a diffusive scaling in \cite{L.S2012}. An arbitrary order version, based on Runge-Kutta methods, has been proposed recently in \cite{L.L.S2016}, where it was also analysed for kinetic equations in the advection-diffusion limit. In \cite{L.M.S2017}, the scheme was used to construct a explicit, flexible, arbitrary order method for general non-linear hyperbolic conservation laws, based on their relaxation to a kinetic equation. Alternative approaches to obtain higher-order \ac{PI} schemes have been proposed in \cite{L.G2007,R.G.K2004}.
These methods align with recent research efforts on numerical methods for multiscale simulation \cite{E.L.R+2007,G.H.K+2003}.

For problems exhibiting more than two time scales, \ac{TPI} was proposed in \cite{G.K2003a}. In these methods, the projective idea is applied recursively. Starting from an inner integrator at the fastest time scale, a \ac{PI} method is constructed with a time step that corresponds to the second-fastest scale. This \ac{PI} method is then considered as the inner integrator of yet another \ac{PI} method at a coarser level. By repeating as required, \ac{TPI} methods construct a hierarchy of projective levels, each using the previous one as its inner step. This idea was explored for linear kinetic equations in \cite{M.S2018}. These methods turn out to have a computational cost that is essentially independent of the stiffness of the collision operator. This property was used with great success in the series of papers \cite{M.R.S2017,M.R.S2019} to develop \ac{PI} and \ac{TPI} methods for the full non-linear \acf{BGK} and Boltzmann equations of single-species rarefied gas dynamics.

\ac{PI} methods are not \acf{AP} methods as such, because the schemes cannot be evaluated explicitly at $\varepsilon=0$ to obtain a classical numerical scheme for the limiting equation. Nevertheless, \ac{PI} and \ac{TPI} methods share important features with \ac{AP} methods. In particular, their computational cost does (in many cases) not depend on the stiffness of the problem. To be precise, it was shown in \cite{M.S2018} for linear kinetic equations that the number of inner time steps at each level of the telescopic hierarchy is independent of the small-scale parameter $\varepsilon$, as is the step size of the outermost integrator. The only parameter in the method which may depend on $\varepsilon$ is the {number} of levels in the telescopic hierarchy. For systems in which the spectrum of the collision operator falls apart into a set of clearly separate clusters (each corresponding to a specific time scale), the number of levels equals the number of spectral clusters. In this regime, the computational cost is completely independent from $\varepsilon$. When the collision operator comprises a continuum of time scales, the number of \ac{TPI} levels increases logarithmically with $\varepsilon$.

The linearisation of the multispecies Boltzmann operator around its equilibria has very similar properties to the classical single species Boltzmann operator, as shown in \cite{D.J.M+2016, B.D2016}. Among others, its spectrum is well separated, with slow modes close to the origin, fast modes at a distance of order $\varepsilon^{-1}$ left of the imaginary axis, and an essential part even farther away from the axis \cite{E.P1975,Rey2013}; this information was the basis for the development of the \ac{TPI} method in \cite{M.R.S2019}. The spectral gap estimates in the multispecies setting are also very similar to the single species case. In this paper, we will exploit this structure to develop a \ac{TPI} method for the multispecies Boltzmann equation.

The rest of the work is organised as follows. In \cref{sec:kinEqMult} we recall the elements of kinetic theory before presenting the models under consideration: the multispecies Boltzmann and \ac{BGK} equations. \Cref{sec:schemes} develops the numerical schemes which we shall use to study these models, with emphasis on the novelty of this paper: the use of \acl{PI} and \acl{TPI} to construct uniformly accurate schemes for the multispecies \ac{BGK} model. In Section \ref{sec:spectra} we compute the approximate spectrum of the numerical scheme and discuss the strategy for the choice of the parameters of the numerical methods. Finally, \cref{sec:experiments} presents a variety of numerical experiments, both in $1 + 1$ and $2 + 2$ dimensions of the phase space, which demonstrate how the schemes can handle extreme mass ratios, fluid instabilities, and other complex phenomena.
\section{On Kinetic Equations and Multiple Species Models}
\label{sec:kinEqMult}

The cornerstone of kinetic theory is the Boltzmann equation, which describes the evolution in time of the distribution of the particles of a rarefied gas. Each particle is subject to ballistic motion, travelling in a straight line at a given velocity, which may only change when it collides with another particle. These collisions are assumed to take place instantaneously and to conserve momentum and kinetic energy. Moreover, the gas is assumed to be dilute enough so that collisions between three or more particles can be neglected.

The \textit{distribution} of particles, $f\e (t,\bx,\bv)$, is a non-negative function which describes for every time $t$ the likelihood of finding a particle \revisionTwo{in the infinitesimal configuration $\prt{\dbx,\dbv}$.} In this work, $f\e$ is understood in the \textit{number} sense, meaning it is not a probability measure. To be precise\footnotemark,
$
	\iint f\e(t,\bx,\bv) \dbx\dbv = N
$
for all times $t\geq 0$, where $N$ is the number of particles in the gas. All particles are assumed identical, each with mass $m > 0$.
\footnotetext{All the integrals in the text are taken over the full position domain, the full velocity domain, or the entire domain of \Cref{eq:kinetic}, unless otherwise stated:
	\begin{align}
		\int \dbx \equiv \int_\Omega \dbx,\quad
		\int \dbv \equiv \int_\Rv \dbv,\quad
		\iint \dbx\dbv \equiv \iint_{\Omega \times \Rv} \dbx \dbv.
	\end{align}
}

The distribution $f\e$ evolves according to the kinetic equation
\begin{align}\label{eq:kinetic}
	\begin{system}
		\pt f\e + \bv\cdot \grad_\bx f\e = \frac{1}{\eps} \Q[f\e],
		\quad \bx\in\Omega\subseteq\Rx,
		\, \bv\in\Rv,
		\, t>0,
		\\
		f\e(0,\bx,\bv) = f_0(\bx,\bv),
	\end{system}
\end{align}
\revisionTwo{for some given Lipschitz domain $\Omega$, which has to be supplemented with boundary conditions. The most simple physically relevant boundary conditions are the so-called Maxwell boundary conditions, where a fraction $1-\alpha \in [0,1]$ of the particles is {specularly} reflected by a wall, whereas the remaining fraction $\alpha$ is {thermalized}, leaving the wall (moving at a velocity $u_w$) at a temperature $T_{w} > 0$. Such conditions can be written as follows: given the unit outer normal ${n}_\bx$ for all $\bx \in \partial \Omega$, for an {incoming} velocity $\bv$ (namely $\bv \cdot {n}_\bx \leq 0$), we set
	$
		f(t,\bx,\bv) = \alpha \,\mathcal{R}f(t,\bx,\bv) \,+\, (1-\alpha) \,\mathcal{M}f(t,\bx,\bv),
	$
	with
	\[
		\left\{ \begin{aligned}
			 & \mathcal{R}f(t,\bx,\bv) \,\,=\,\, f(t,\bx, \bv\,-\, 2 ({n}_\bx \cdot \bv)\, {n}_\bx),
			\\
			 & \mathcal{M}f(t,\bx,\bv) \,\,=\,\, \frac{\mu(t,\bx)}{(2 \pi \, T_w)^{\dbv/2}} \exp\left (-\frac{|\bv - u_w|^2}{2\,T_{w}}\right ).
		\end{aligned}\right.
	\]
	In this last expression, $\mu$ insures global mass conservation.
	Nevertheless, for the sake of simplicity, we shall only consider either periodic or influx boundary conditions.}
This equation couples a \textit{transport operator}, which models the displacement of particles, with a \textit{collision operator} $\Q$, which describes the changes in momentum due to collisions. In this work, $\dimx=\dimv=1$ or $2$.

The frequency of collisions is governed by the \revisionTwo{dimensionless \textit{Knudsen number} $\varepsilon$, defined as the ratio between the mean free path of particles and the length scale of observation; broadly, }a measure of the average distance that a particle can travel before colliding with another one. The value of $\varepsilon$ distinguishes the \textit{kinetic regime}, where collisions are rare, from the \textit{hydrodynamic regime}, where collisions dominate the dynamics; the former corresponds to larger values of $\varepsilon$, and the latter, to smaller values.

The seminal Boltzmann collision operator $\Q$
is a quadratic operator local in $(t,\bx)$. The time and position act only as parameters in $\Q$ and therefore will be omitted in its description. It is given, for $\dimv \geq 2$, by
\begin{equation}
	\label{eq:boltzmannOperator}
	\Q [f](\bv) = \iint_{\Rv \times\mathbb{S}^{\dimv-1}} B(|\bv-\bv_*|,\cos \theta) \,
	\left( f' f'_* - f f_* \right) \dbv_* \dd \sigma,
\end{equation}
where we have used the shorthand $f = f(\bv)$, $f_* = f(\bv_*)$,
$f ^{'} = f(\bv')$, $f_* ^{'} = f(\bv_* ^{'})$. The velocities of the
colliding pairs, $(\bv,\bv_*)$ and $(\bv',\bv'_*)$, can be parametrized as
\begin{equation}
	\bv' = \frac{\bv+\bv_*}{2} + \frac{|\bv-\bv_*|}{2} \sigma, \qquad
	\bv'_* = \frac{\bv+\bv^*}{2} - \frac{|\bv-\bv_*|}{2} \sigma.
\end{equation}
The collision kernel $B$ is a non-negative function which, by physical arguments of invariance, may only depend on $|\bv-\bv_*|$ and
$\cos \theta = {\hat g} \cdot \sigma$ (where ${\hat g} =
	(\bv-\bv_*)/|\bv-\bv_*|$).
It characterises the
details of the binary interactions, and has the form
\begin{equation}
	\label{defVHSKernel}
	B(|\bv-\bv_*|,\cos\theta)=|\bv-\bv_*|\,\Phi(|\bv-\bv_*|,\cos\theta).
\end{equation}
The {scattering cross-section} $\Phi$, in the case of inverse $k$\th power forces between particles, can be written as
\[ \Phi(\vert \bv - \bv_{\ast} \vert,
	\cos\theta) = b_{\alpha}(\cos\theta) \, \vert \bv - \bv_{\ast}
	\vert^{\alpha-1},
\]
with $\alpha=(k-5)/(k-1)$.
The special situation {$k=5$} gives the so-called {Maxwell pseudo-molecules model} with
$ B(|\bv-\bv_*|,\cos \theta)=
	b_{0}(\cos\theta). $
For the Maxwell case the collision kernel is independent of the relative velocity. For numerical purposes, a
widely used model is the {variable hard sphere} (VHS)
model introduced by Bird~\cite{Bird1994}. The model corresponds to {$b_{\alpha}(\cos\theta)=C_\alpha$,}
where {$C_\alpha$} is a positive constant, and hence
$ \Phi(\vert \bv - \bv_{\ast} \vert, \cos\theta) =
	C_{\alpha} \vert \bv - \bv_{\ast} \vert^{\alpha-1}. $
For further details on the physical background and derivation of
the Boltzmann equation, we refer to~\cite{C.I.P1994, Villani2002}

\Cref{eq:kinetic} is rich with structural properties. In terms of conservation, its weak form reveals that any conserved quantity of the dynamics corresponds to an associated \textit{collision invariant} of $\Q$: a function $\varphi(\bv)$ such that
\begin{equation}
	\int \Q[f] (\bv) \, \varphi(\bv) \dbv = 0
\end{equation}
for any function $f(\bv)$. For the Boltzmann operator, the space of collision invariants is given by $\operatorname{span} \set{1, \bv, \abs{\bv}^2}$ \cite{Cercignani1990}. In consequence, the \textit{mass}, \textit{momentum}, and \textit{kinetic energy} of the distribution, respectively
$\iint m f\e(t,\bx,\bv) \dbx\dbv$,
$\iint m \bv f\e(t,\bx,\bv) \dbx\dbv$,
and $\iint \frac{m}{2} \abs{\bv}^2 f\e(t,\bx,\bv) \dbx\dbv$,
are conserved.

In terms of dynamics, the relaxation properties of \cref{eq:kinetic} are well-understood. The celebrated H-theorem characterises the dissipation of the \textit{Boltzmann entropy}\footnotemark,
\begin{align}
	\curlyH = \int f\e(\bv) \log(f\e(\bv)) \dbv.
\end{align}
At the local level, the theorem states
\begin{align}\label{eq:HTheorem}
	\curlyS = \int \Q[f\e](\bv) \log(f\e(\bv)) \dbv \leq 0,
\end{align}
where $\pder{\curlyH}{t} + \divx \curlyJ = \frac{\curlyS}{\varepsilon}$ and $\curlyJ = \int \bv f\e(\bv) \log(f\e(\bv)) \dbv$, see, for instance, \cite{Cercignani1988,C.I.P1994}. Furthermore, the equality in \eqref{eq:HTheorem} is only achieved for states in the kernel of $\Q$. Such equilibrium states are always Maxwellian distributions:
\begin{align}
	f\e = \M^{n,\barbv,T}(t,\bx,\bv) \coloneqq
	n(t,\bx)\, \prt*{\frac{m}{2\pi T}}^{\dimv/2} \exp\prt*{-\frac{m(\barbv-\bv)^2}{2 T}}
\end{align}
for particles of mass $m$. The \textit{number density} $n$, the \textit{mass density} $\rho$, the \textit{average velocity} $\barbv$, the \textit{temperature} $T$, and the \textit{pressure} $\curlyP$ are (local) moments computed from the distribution which depend only on time and position; they are given by
\begin{align}\label{eq:moments}
	n = \int f\e \dbv,\quad
	\rho = m n,\quad
	\rho\barbv = \int m\bv f\e \dbv,\quad
	\frac{\dimv}{2}nT = \int \frac{m}{2} \abs{\barbv - \bv}^2 f\e \dbv,
\end{align}
where $\rho = m n$ and $\curlyP=nT$.
Note, in particular, that the Maxwellian $\M^{n,\barbv,T}$ has moments $n,\barbv,T$.
\footnotetext{
	\revisionTwo{As is usual, we choose units such that Boltzmann's constant equals one, $k_B=1$, throughout the text.}
}

\subsection{The BGK Model}

The collision operator \eqref{eq:boltzmannOperator} has a complicated structure from both the analytical and the numerical perspectives. Because of this, many works have proposed simpler operators which attempt to capture some or all of the structural properties of $\Q$. A ubiquitous simplification is the model of \acl{BGK} \cite{B.G.K1954}; since the overall effect of the Boltzmann operator is to drive $f\e$ towards the corresponding Maxwellian, they propose a reduced operator which makes the relaxation explicit:
\begin{align}\label{eq:BGK}
	\QBGK[f\e] = \nu(\curlyM[f\e] - f\e),
\end{align}
where $\curlyM[f\e]$ is the Maxwellian whose moments $n,\barbv,T$ are the ones of $f\e$, and where $\nu(t, \bx)$ is a positive collision rate to be determined. This simplification can be derived from $\Q$ by assuming that the distribution $f\e$ is already close to equilibrium. The \ac{BGK} equation is thus defined as the evolution law \eqref{eq:kinetic} together with the operator \eqref{eq:BGK}.

Despite its apparent simplicity, the \ac{BGK} model captures many structural aspects of the Boltzmann collision operator. Because $\curlyM[f\e]$ and $f\e$ share the first three moments, the invariants of $\Q$ are also invariants of $\QBGK$ whenever the collision rate does not depend on $\bv$. Therefore, the new dynamics conserve the mass, momentum, and kinetic energy of $f\e$ as well. The dynamical structure of $\Q$ is also preserved by the \ac{BGK} model because the H-theorem still holds:
\begin{align}\label{eq:HTheoremBGK}
	\curlyS = \int \QBGK[f\e](\bv) \log(f\e(\bv)) \dbv \leq 0,
\end{align}
see, for instance, \cite{Struchtrup2005}. Once again, the equality is only achieved by the equilibrium states, which are the Maxwellians by construction.

\subsection{The Multispecies Boltzmann Equation}
\label{sec:multispeciesBoltz}

A limitation of the Boltzmann equation \eqref{eq:kinetic} is the assumption that all gas particles are identical. However, it is possible to extend this equation to the case of a mixture of gases by the physical arguments used in the single species case \cite{Briant2015}. As stated in the introduction, the derivation of the multispecies Boltzmann equation is mostly formal, unlike that of the single species case, which has been clearly established in a variety of works such as \cite{G.S.T2013}.

We will consider a mixture of $P$ species, each described by a distribution\footnotemark~$f\p (t,\bx,\bv)$, a time-dependent non-negative function as before. Each distribution is understood in the number sense:
\begin{align}
	\iint f\p(t,\bx,\bv) \dbx\dbv = N\p
\end{align}
for all times $t\geq 0$, where $N\p$ is the number of particles in the $p\th$ species. All particles of that species are assumed identical, each with mass $m\p$.
\footnotetext{The $\varepsilon$ notation is dropped in the interest of simplicity. Henceforth, $f$ stands for $f\e$.
}

The distributions $f\p$ evolve according to the so-called multispecies Boltzmann equation, given by
\begin{align}\label{eq:kineticMultispecies}
	\begin{system}
		\pt f\p + \bv\cdot \grad_\bx f\p = \frac{1}{\eps} \Q\p[\bbf],
		\quad \bx\in\Omega\subseteq\Rx,\, \bv\in\Rv,\, t>0,\\
		\Q\p[\bbf] = \sum_{q=1}^{P} \Q\pq[f\p,f\q],
		\\
		f\p(0,\bx,\bv) = f_{p,\,0}(\bx,\bv),
	\end{system}
\end{align}
posed with appropriate boundary conditions. The changes of momentum of the $p\th$ species are now governed by a sum of collision operators, one for each species in the gas.
These multispecies collision operators $\Q\pq$ are similar to the classical Boltzmann operator \eqref{eq:boltzmannOperator}; namely, they are given by
\begin{equation}
	\label{eq:boltzmannOperatorMultispecies}
	\Q\pq[f\p,f\q] = \iint_{\Rv \times\mathbb{S}^{\dimv-1}} B\pq(|\bv-\bv_*|,\cos \theta)
	\left[ f\p\left(\bv'\right) f\q\left( \bv'_*\right) - f\p(\bv)\, f\q(\bv_*) \right] \dbv_* \dd \sigma,
\end{equation}
where inter-species collisions are given by
\begin{align}\label{eq:collisionMultiSpec}
	\begin{split}
		& \bv' \,=\, \frac{1}{m\p + m\q}\left (m\p \bv +m\q \bv_* +m\q |\bv-\bv_*| \sigma\right ), \\
		& \bv'_* \,=\, \frac{1}{m\p + m\q}\left (m\p \bv +m\q \bv_* -m\p |\bv-\bv_*| \sigma\right ).
	\end{split}
\end{align}
Note that in the intra-species collision case $p = q$,
\eqref{eq:boltzmannOperatorMultispecies} corresponds to the classical Boltzmann operator \eqref{eq:boltzmannOperator}. We shall also assume for the sake of simplicity that the inter-species collision kernel $B\pq$ is independent on $p$ and $q$, and equal to the variable hard spheres kernel \eqref{defVHSKernel}.

The conservation properties of the Boltzmann equation persist in the multispecies case at the level of the mixture. The mass of each species,
$
	\iint m\p f\p(t,\bx,\bv) \dbx\dbv,
$
is conserved because $1$ remains an invariant of each operator:
\begin{align}
	\int \Q\pq[f,g] (\bv) \dbv = 0
\end{align}
for any functions $f(\bv)$ and $g(\bv)$. However, the functions $\bv$ and $\abs{\bv}^2$ are only invariants in the sense that
\begin{align}\label{eq:balanceMomentum}
	\int \bv \Q\pq[f,g] (\bv) \dbv = - \int \bv \Q\qp[g,f] (\bv) \dbv
\end{align}
and
\begin{align}\label{eq:balanceEnergy}
	\int \abs{\bv}^2 \Q\pq[f,g] (\bv) \dbv = - \int \abs{\bv}^2 \Q\qp[g,f] (\bv) \dbv.
\end{align}
As such, it is only the \textit{total momentum} and the \textit{total kinetic energy} of the gas, respectively
\begin{align}
	\sum\limits_{p=1}^{P} \iint m\p \bv f\p(t,\bx,\bv) \dbx\dbv,
	\text{ and }\,
	\sum\limits_{p=1}^{P} \iint \frac{m\p}{2} \abs{\bv}^2 f\p(t,\bx,\bv) \dbx\dbv,
\end{align}
that are conserved.

The H-theorem also applies to the multispecies case, for the \textit{total Boltzmann entropy} and total dissipation:
\begin{align}\label{eq:HTheoremMultispecies}
	\curlyH = \sum\limits_{p=1}^{P} \int f\p(\bv) \log(f\p(\bv)) \dbv,\quad
	\curlyS = \sum\limits_{p=1}^{P} \int \Q\p[f\p](\bv) \log(f\p(\bv)) \dbv \leq 0,
\end{align}
see \cite{C.C1990}. The equilibrium distributions are all Maxwellians with common average velocity, $\barbv\eq$, and temperature, $T\eq$:
\begin{align}
	\label{eq:equilibriumMaxwellianfp}
	f\p = \M^{n,\barbv\eq,T\eq}(t,\bx,\bv) = n\p(t,\bx)\, \prt*{\frac{m\p}{2\pi T\eq}}^{\dimv/2} \exp\prt*{-\frac{m\p(\barbv\eq-\bv)^2}{2 T\eq}}.
\end{align}
Note $n\p$ and the rest of individual moments are defined just as in \cref{eq:moments}:
\begin{align}\label{eq:momentsMultispecies}
	\begin{split}
		n\p = \int f\p \dbv,\quad
		\rho\p\barbv\p = \int m\p\bv f\p \dbv,\quad
		\frac{\dimv}{2}n\p T\p = \int \frac{m\p}{2} (\barbv\p - \bv)^2 f\p \dbv,\end{split}
\end{align}
where $\rho\p = m\p n\p$ and $\curlyP_p=n_pT_p$.
The total moments for the mixture are given by
$n = \sum n\p$,
$\rho = \sum \rho\p$,
$\rho\barbv = \sum \rho\p\barbv\p$,
$\curlyP = \sum \curlyP\p$, and
\begin{align}\label{eq:momentsMultispeciesOverall}
	\begin{split}
		\frac{\dimv}{2}n T = \sum\limits_{p=1}^{P} \int \frac{m\p}{2} (\barbv - \bv)^2 f\p \dbv.
	\end{split}
\end{align}

\subsection{A Multispecies BGK Model}

Although the multispecies Boltzmann equation described in \cref{sec:multispeciesBoltz} shares many features with the single species Boltzmann equation, the fact that conservations are only global in space and all species, e.g. \eqref{eq:balanceEnergy}, makes its simplification to a relaxation operator such as the \ac{BGK} model \eqref{eq:BGK} much harder. Many different approaches exists to derive such a multispecies BGK operator, but very few are satisfactory in regards to the macroscopic properties of the resulting relaxation operator. The interested reader can find more details about these various approaches in the recent papers \cite{B.I.P2019,B.P2020}.

This work will follow along the lines of \cite{H.H.M2017}. The formulation of the multispecies \ac{BGK} operator requires the definition of the mixture Maxwellians:
\begin{align}\label{eq:mixtureMaxwellian}
	\curlyM\pq[f\p, f\q](t,\bx,\bv) = n\p(t,\bx) \prt*{\frac{m\p}{2\pi T\pq}}^{{\dimv}/{2}} \exp\prt*{-\frac{m\p(\barbv\pq-\bv)^2}{2 T\pq}},
\end{align}
where $\barbv\pq$ and $T\pq$ are mixture moments to be determined. It is assumed that the interaction of the $p\th$ and the $q\th$ species drives $f\p$ towards $\curlyM\pq$, as well as $f\q$ towards $\curlyM\qp$; thus we require $\barbv\pq=\barbv\qp$, $T\pq=T\qp$. Assuming the distributions are close to equilibrium permits the simplification of $\Q\pq$ to obtain
\begin{align}\label{eq:BGKMultispecies}
	\QBGK\pq[f\p,f\q] = \nu\pq(\curlyM\pq[f\p, f\q] - f\p),
\end{align}
where $\nu\pq$ are positive collision rates\footnotemark. The multispecies \ac{BGK} model is thus defined as the evolution law \eqref{eq:kineticMultispecies} together with the operator \eqref{eq:BGKMultispecies}.
\footnotetext{
	\revisionTwo{
		From the derivation of the multispecies \ac{BGK} model,
		\begin{equation}
			\nu\pq(\bv) = \iint_{\Rv \times\mathbb{S}^{\dimv-1}} B\pq(|\bv-\bv_*|,\cos \theta)
			\curlyM\qp(\bv)
			\dbv_* \dd \sigma,
		\end{equation}
		see \cite{H.H.M2017}. However, as done in the reference, we assume the collision frequencies are given constants, $\nu\pq\equiv 1$.
	}
}

In order to recover the structural properties of the Boltzmann equation, the mixture moments have to be chosen appropriately. Considering the interaction of the $p\th$ species with itself, we deduce $\barbv\pp=\barbv\p$ and $T\pp=T\p$. To define $\barbv\pq$, it is enough to impose the conservation of momentum in the form of \cref{eq:balanceMomentum}. This readily yields
\begin{align}\label{eq:mixtureVelocity}
	\barbv\pq = \frac{\rho\p \nu\pq \barbv\p + \rho\q \nu\qp \barbv\q}{\rho\p \nu\pq + \rho\q \nu\qp}.
\end{align}
This expression is reminiscent to the microscopic parametrisation of the postcollisional velocities as a function of the precollisional ones \eqref{eq:collisionMultiSpec}.
Similarly, requiring \cref{eq:balanceEnergy} as an expression of conservation of kinetic energy produces
\begin{align}\label{eq:mixtureTemperature}
	T\pq =
	\frac{n\p \nu\pq T\p + n\q \nu\qp T\q}{n\p \nu\pq + n\q \nu\qp} +
	\frac{\rho\p \nu\pq (\barbv\p^2 - \barbv\pq^2) + \rho\q \nu\qp (\barbv\q^2 - \barbv\qp^2)}{\dimv(n\p \nu\pq + n\q \nu\qp)}.
\end{align}
This expression guarantees $T\pq\geq 0$, see \cite[Appendix 9]{H.H.M2017}. These choices ensure, by construction, the conservation of total momentum and total kinetic energy of the gas.

This model also captures the dissipative structure of the multispecies Boltzmann equations. It is shown in \cite{H.H.M2017} that the collision operator \eqref{eq:BGKMultispecies} verifies the H-theorem:
\begin{align}\label{eq:HTheoremMultispeciesBGK}
	\curlyS = \sum\limits_{p=1}^{P} \int \QBGK\p[f\p](\bv) \log(f\p(\bv)) \dbv \leq 0.
\end{align}

\subsection{Hydrodynamic Limits}\label{sec:hydrolimits}

We will briefly present here the zeroth order hydrodynamic limit of the multispecies BGK model \eqref{eq:BGKMultispecies}, following along the lines of the recent paper \cite{B.B.G+2019}. Taking the limit $\varepsilon \to 0$ in that kinetic equation will formally project the distributions $f\p$ towards the equilibrium distributions \eqref{eq:equilibriumMaxwellianfp} whose moments solve a system of conservation laws. This system is obtained by multiplying the kinetic multispecies equation \eqref{eq:kineticMultispecies} by the powers of $\bv$, integrating in velocity, and summing the equations on the individual moments in order to have only $\dimv$ equations for the total momentum.

Under the Maxwellian closure, we obtain the following system of $\dimv+P+1$ conservation law, the so-called multispecies Euler limit of \eqref{eq:kineticMultispecies}:
\begin{equation}
	\label{eq:MultispeciesEuler}
	\begin{system}
		\partial_t n\p + \grad_\bx \cdot \left ( n\p \barbv \right ) = 0, \quad p\in \{1, \cdots, P\}, \\[.5em]
		\partial_t \left ( \rho \, \barbv \right ) + \grad_\bx \cdot \left ( \rho \, \barbv \otimes \barbv \right ) + \grad_\bx \left( n \, T \right) = \mathbf{0}, \\[.5em]
		\partial_t \left( \frac 12 \rho \, |\barbv|^2 +\frac \dimv{2} n \, T \right ) + \grad_\bx \cdot \left [\left ( \frac 12 \rho \, |\barbv|^2 +\frac {\dimv+2}{2} n \, T \right ) \barbv \right] = 0.
	\end{system}
\end{equation}
The interested reader can consult \cite{B.B.G+2019} for the terms corresponding to the next order in the expansion, the multispecies Navier-Stokes limit of \eqref{eq:kineticMultispecies}. \section{Numerical Schemes}\label{sec:schemes}

This section is devoted to the development of efficient numerical methods for the multispecies \ac{BGK} model, given by \cref{eq:kineticMultispecies,eq:BGKMultispecies}. The numerical simulation of the model is hindered by the presence of two radically different scales. The transport term classically imposes a CFL condition on the mesh size: \revisionOne{$\Dt = \mathcal{O}\prt{\Dx}$} is needed for stability. However, in the regimes where $\varepsilon\ll 1$, the collision term dominates the dynamics, further restricting the time step to \revisionTwo{$\Dt = \mathcal{O}\prt{\varepsilon}$}. This requirement is particularly pernicious in view of the large cost per step (due to the computation of non-local terms), and it would render a na\"{i}ve scheme inefficient. The development of schemes for kinetic and hyperbolic equations which can handle multiple scales without such penalties (sometimes called \textit{asymptotic-preserving} schemes) is a very active area of research, see the recent reviews \cite{Jin2012,D.P2014,H.J.L2017}.

In order to construct an efficient numerical scheme for \Cref{eq:kineticMultispecies}, we resort to \textit{projective} and \textit{telescopic projective integration}. First introduced in \cite{G.K2003}, projective integration provides an efficient framework for the numerical solution of differential equations with two distinct timescales, characterised by two separate clusters of eigenvalues. Telescopic projective integration was introduced by \cite{G.K2003a} in order to extend the method to problems with many scales, or without a clear spectral separation. These methods are not asymptotic preserving in the classical sense \cite{Jin1999} because they cannot be evaluated at $\varepsilon = 0$ to recover a scheme for the asymptotic equation; nevertheless, they are able to solve stiff problems efficiently, and have been successfully adapted to some kinetic equations: radiative transfer \cite{L.S2012}, the linearised \ac{BGK} equation with multiple relaxation times \cite{M.S2018}, and both the full \ac{BGK} model and the Boltzmann equation \cite{M.R.S2017,M.R.S2019}. Higher-order approaches have also been developed \cite{L.L.S2016,L.M.S2017}.

In the remainder of this section we describe the numerical method. First we introduce the phase space discretisation, using a discrete velocity method for the collision operator, coupled with a finite-volume scheme for the transport. Then we present the projective and telescopic schemes.

\subsection{Phase Space Discretisation}

\paragraph{Velocity}

The discretisation of \Cref{eq:kineticMultispecies} can be conducted in a number of ways; see \cite{D.P2014} for a survey. Here we shall employ a discrete velocity approach as in \cite{R.S1994,Buet1996,Mieussens2000}, where the velocity space is approximated by a finite grid. The choice of grid is non-trivial, and it can affect the entropic properties of the numerical scheme \cite{A.P.P2008}; for simplicity, we choose a Cartesian grid.

The velocity space, $\R^\dimv$, is first restricted to a bounded domain $\prt{-L_v,L_v}^\dimv$. Each dimension is discretised into $2N_v+1$ points, separated by a distance $\Dv = 2L_v/(2N_v+1)$; the $j\th$ point is given by $v\j=j\Dv$, where $j\in\curlyN_v\coloneqq\set{-N_v,\cdots,N_v}$. The discretised velocity space becomes $\curlyV = \set{\bv_\bj = \prt{v_{j_1},\cdots,v_{j_\dimv}}}$, where $\bj=\prt{j_1,\cdots,j_\dimv}\in \curlyN_v^\dimv$ is a multi-index. \Cref{eq:kineticMultispecies} thus becomes a system of transport equations:
\begin{align}\label{eq:kineticMultispeciesDV}
	\begin{system}
		\pt f\pbj + \bv_{\bj}\cdot \grad_\bx f\pbj = \frac{1}{\eps} \Q\pbj[f\p],
		\quad \bx\in\Omega\subseteq\Rx,\, \bj\in\curlyN_v^\dimv,\, t>0,\\
		\Q\pbj[f\p] = \sum_{q=1}^{P} \Q\pqbj[f\p,f\q],
	\end{system}
\end{align}
where $f\pbj$ approximates $f\p(\bv_\bj)$, and $\Q\pqbj[f\p,f\q]$ approximates $\Q\pq[f\p,f\q](\bv_\bj)$. The initial datum $f\pbj(0,\bx)$ is prescribed by evaluating the continuous datum, $f_{p,\,0}$, at $\bv_\bj$. Note that the evolution of $f\pbj$ depends on $\Q\pbj$, which itself depends on the entire distribution of every species, coupling all the equations in the system.

The evaluation of the collision operator will require the approximation of the moments of the distribution based on the discretised velocity grid.
The semi-discrete individual moments are\footnotemark
\begin{align}\label{eq:momentsMultispeciesDiscrete}
	n^h\p = \sum f\pbj \Dbv,\quad
	\rho^h\p\barbv^h\p = \sum m\p\bv_\bj f\pbj \Dv,\quad
	\frac{\dimv}{2}n^h\p T^h\p = \sum \frac{m\p}{2} (\barbv^h\p - \bv_\bj)^2 f\pbj \Dv,
\end{align}
where $\Dbv=\Dv^\dimv$, and $\rho^h\p = m\p n^h\p$ as before. The approximate mixture moments are computed as in \cref{eq:mixtureVelocity,eq:mixtureTemperature}. The collision operator is, therefore,
\begin{align}\label{eq:BGKMultispeciesDV}
	\QBGK\pqbj[f\p,f\q] = \nu\pq(\curlyM\pqbj[f\p, f\q] - f\pbj),
\end{align}
where $\curlyM\pqbj[f\p, f\q] = \curlyM\pq[f\p, f\q](\bv_\bj)$, and $\curlyM\pq[f\p, f\q]$ is the Maxwellian \eqref{eq:mixtureMaxwellian} generated by the approximate moments.
\footnotetext{
	As with the integrals, the sums in this text are taken over the full velocity domain, unless otherwise stated:
	\begin{align}
		\sum \Dbv \equiv \sum_{\bj\in \curlyN_v^\dimv} \Dbv,
	\end{align}
	\revisionTwo{where $\curlyN_v\coloneqq\set{-N_v,\cdots,N_v}$, see the discussion around \cref{eq:kineticMultispeciesDV}.}
}

\paragraph{Space}

The spatial discretisation of \Cref{eq:kineticMultispeciesDV} will be guided by the treatment of the transport term alone. Linear transport terms such as $\bv_{\bj}\cdot \grad_\bx f\pbj$ can be discretised in a number of ways, including semi-Lagrangian methods \cite{C.R.S2009,C.M.S2010,D.L.N+2018} and \revisionOne{finite-volume methods \cite{LeVeque1990,Toro1999,E.G.H2000,LeVeque2002}}. Here we shall prescribe a first-order finite-volume scheme, but this can be replaced with a higher-resolution method (a WENO scheme \cite{L.O.C1994}, for instance) independently of the choices for the discretisations of velocity and time.

The spatial domain, $\Omega$, is assumed to be rectangular, of the form $\prt{0,L_{x,1}} \times \cdots \times \prt{0,L_{x,\dimx}}$. Each dimension $d$ is discretised into $N_{x,d}$ volumes of size $\Dx_{d} = L_{x,d}/N_{x,d}$: every one, a cell $C\i = \prt{x\imh,x\ih}$ with centre $x\i$, where $x\i=i\Dx$, for $i\in\curlyN_x\coloneqq\set{1,\cdots,N_x}$. The discretised spatial domain becomes $\curlyX=\set{\bC_\bi=C_{i_1}\times\cdots\times C_{i_\dimx}}$, where $\bi=\prt{i_1,\cdots,i_\dimx}\in \curlyN_x^\dimx$ is a multi-index. \Cref{eq:kineticMultispeciesDV} thus becomes a system of equations:
\begin{align}\label{eq:kineticMultispeciesDVDX}
	\begin{system}
		\pt f\pbibj
		+ \Phi\pbibj[f\p]
		= \frac{1}{\eps} \Q\pbibj[f\p],
		\quad \bi\in\curlyN_x^\dimx,\, \bj\in\curlyN_v^\dimv,\, t>0,\\
		\Q\pbibj[f\p] = \sum_{q=1}^{P} \Q\pqbibj[f\p,f\q],
	\end{system}
\end{align}
where $f\pbibj$ is a finite-volume approximation of $f\pbj$,
$
	f\pbibj(t) \simeq \frac{1}{\abs{\bC_\bi}} \int_{\bC_\bi} f\pbj(t,\bx) \dbx,
$
with equality for the initial datum.

The term $\Phi\pbibj$, expresses the net flux across the boundary of $\bC_\bi$. In the one-dimensional setting, $\dimx=1$, it is written in conservation form as
\begin{align}
	\Phi\pij[f\p] = \frac{F\pihj-F\pimhj}{\Dx},
\end{align}
where
\begin{align}
	F\pihj = f\pij \pos{v\j} + f\pipj \neg{v\j}
\end{align}
is the upwind flux discretisation. \revisionOne{
	Respectively, $\pos{v\j}$ and $\neg{v\j}$ are the positive and negative parts of $v\j$; $\pos{s}=\neg{-s}=\max\set{s,0}$ for any $s\in\R$.
} The generalisation to higher dimensions is immediate, \revisionOne{see any of \cite{LeVeque1990,Toro1999,E.G.H2000,LeVeque2002}}.

The discretisation in space replaces the collision operator \eqref{eq:BGKMultispeciesDV} with
\begin{align}\label{eq:BGKMultispeciesDVDX}
	\QBGK\pqbibj[f\p,f\q] = \nu\pq(\curlyM\pqbibj[f\p, f\q] - f\pbibj),
\end{align}
where $\curlyM\pqbibj[f\p, f\q]$ is simply the Maxwellian computed from the discretisations of $f\p, f\q$ in velocity and space. The moments given in \eqref{eq:momentsMultispeciesDiscrete} are easily discretised;
for example,
\begin{align}
	n\p(t,\bx) = \sum f\pbj(t,\bx) \Dbv
	\quad\text{becomes}\quad
	n\pbi(t) = \sum f\pbibj(t) \Dbv.
\end{align}

\subsection{Projective Integration}

Following the phase space discretisation, the multispecies \ac{BGK} equation leads to the semi-discrete system \eqref{eq:kineticMultispeciesDVDX}. This can be summarised as
\begin{align}\label{eq:semidiscrete}
	\pt f\pbibj = \D\e[f\pbibj],
\end{align}
where $p$, $\bi$, and $\bj$ are, respectively, the species, spatial, and velocity indices, and where $\D\e$ is an operator which comprises both discrete transport and collisions:
\begin{align}
	\D\e[f\pbibj] \coloneqq - \Phi\pbibj[f\p] + \frac{1}{\eps} \Q\pbibj[f\p].
\end{align}

Projective integration essentially consists of two stages: an inner step, where an elementary time-stepping method (the \textit{inner integrator}) is used to solve \eqref{eq:semidiscrete} over a short interval; and a projective step, where the previous information is used to approximate the time derivative of $f$ and extrapolate forward (through the \textit{outer integrator}) over a large time interval.

We define step sizes $\Dt_0$ and $\Dt_1$, respectively for the inner and outer levels. The solution $f\pbibj(t)$ is approximated by $f\pbibj\projnk$ at the time $t=n\Dt_1 + k\Dt_0$. The evolution is computed as follows:

\paragraph{Inner Integrator}

The inner step is computed through the forward Euler method with step size $\Dt_0$:
\begin{align}\label{eq:innerIntegrator}
	f\projnkp = \FE_{\Dt_0}[f\projnk] \coloneqq f\projnk + \Dt_0 \D\e[f\projnk],
\end{align}
where we have omitted the $p$, $\bi$, and $\bj$ indices for simplicity. The inner step size will be chosen as \revisionTwo{$\Dt_0 = \mathcal{O}\prt{\varepsilon}$} to ensure the stability of \eqref{eq:innerIntegrator}, just as it would be in a classical scheme.

While higher-order inner integrators could easily be constructed, it is shown in \cite{M.S2018} that forward Euler is the choice which gives the projective scheme the best stability properties. If high-order accuracy is desired, this can be achieved at the outer level.

\paragraph{Outer Integrator}

The outer step computes the update $f\projnpzero$ from $f\projnzero$, which can be simply be understood as the steps $f\np$ and $f\n$ of a time discretisation $t=n\Dt_1$. The outer step size will be chosen as \revisionOne{$\Dt_1 = \mathcal{O}\prt{\Dx}$}.

In order to compute the update, the method first computes a sequence of $K+1$ inner updates from $f\projnzero$, using the inner integrator. The last two updates, $f\projnK$ and $f\projnKp$, are used to project the solution forward:
\begin{align}\label{eq:outerIntegrator}
	f\np = \PFE_{\Dt_1}[f\n] \coloneqq f\projnKp + (\Dt_1 - (K+1)\Dt_0) \frac{f\projnKp - f\projnK}{\Dt_0};
\end{align}
the relative size of the extrapolation is
$
	M = \frac{\Dt_1}{\Dt_0} - (K+1) .
$

Since the extrapolation is performed with a first-order approximation of the derivative, this is known as the \textit{projective forward Euler method}. Higher order methods can be constructed by using more points to approximate the time derivative, in the vein of Runge-Kutta methods, as discussed in \cite{G.K2003}.

A schematic summary of the method can be found in \cref{alg:projectiveIntegration}. The exact choice of $K$, $\Dt_0$, and $\Dt_1$ will be discussed in \cref{sec:spectra}.
\begin{figure}[ht]
	\centering
	\begin{minipage}{.6\linewidth}
		\begin{algorithm}[H]
			\SetAlgoNoLine
			\KwData{$N$, $K$, $\Dt_1$, $\Dt_0$}
			$f^0 \leftarrow f(0)$ \;
			\For{$0\leq n\leq N-1$}{
				$f\projnzero \leftarrow f\n$ \;
				\For{$0\leq k\leq K$}{
					$f\projnkp \leftarrow \FE_{\Dt_0}[f\projnk]$ \;
				}
				$f\np \leftarrow \PFE_{\Dt_1}[f\n]$ \;
			}
			\caption{Projective Integration}
			\label{alg:projectiveIntegration}
		\end{algorithm}
	\end{minipage}
\end{figure}

\paragraph{Efficiency of the Method}

As discussed in \cite{G.K2003a}, the efficiency of projective integration can be measured by comparing the number of inner steps required to integrate over a time interval directly to that required by the projective method. Assuming a homogeneous cost per step and negligible overhead arising from the projection step (which is reasonable as the projective step does not involve the computation of moments of $f$), projective integration reduces the overall computational cost of a simulation by a factor of
\begin{align}
	S
	= \frac{\Dt_1 \prt*{\Dt_0}^{-1}}{K+1}
	= \frac{M+K+1}{K+1}.
\end{align}

\subsection{Telescopic Integration}

The stability properties of projective integration are well-understood: the method is suited for problems with exactly two timescales which are spectrally separate. There are, however, many multi-scale problems of interest whose spectra lack the sufficient structure: for instance, the linearised and the full \ac{BGK} models, as soon as the collision rate $\nu$ depends on space or on the distribution $f$ \cite{M.S2018,M.R.S2019,K.S2020a}. This is also the case for the multi-species \ac{BGK} model (\ref{eq:kineticMultispecies},~\ref{eq:BGKMultispecies}), as will be shown in \cref{sec:spectra}.

Telescopic integration extends the projective method to include such problems. The method nests several levels: an innermost level, consisting, as before, of an elementary integrator; and a number of outer levels, each performing a projection step based on the iterations of the level below it. We briefly discuss an approach with two projective levels, which will prove sufficient for our needs, though the method can easily be extended.

We define step sizes $\Dt_0$, $\Dt_1$, and $\Dt_2$, respectively for the inner, middle, and outer levels. The solution $f\pbibj(t)$ is approximated by $f\pbibj\projnkl$ at the time $t=n\Dt_2 + k\Dt_1 + l\Dt_0$. The evolution is computed as follows:

The step of the outer integrator is given by
\begin{align}\label{eq:outerTIntegrator}
	f\np = \PFE_{\Dt_2}[f\n] \coloneqq f\telnKp + (\Dt_2 - (K_1+1)\Dt_1) \frac{f\telnKp - f\telnK}{\Dt_1};
\end{align}
the relative size of the extrapolation is
$
	M_1 = \frac{\Dt_2}{\Dt_1} - (K_1+1).
$
In order to evaluate the update, $K_1+1$ steps of the middle integrator have to be performed, to find $f\telnK$ and $f\telnKp$ from $f\n$. Each of those steps is given by
\begin{align}\label{eq:middleTIntegrator}
	f\projnkp = \PFE_{\Dt_1}[f\projnk] \coloneqq f\projnkLp + (\Dt_1 - (K_0+1)\Dt_0) \frac{f\projnkLp - f\projnkL}{\Dt_0},
\end{align}
extrapolations of size
$
	M_0 = \frac{\Dt_1}{\Dt_0} - (K_0+1).
$
For each step to be performed, $K_0+1$ steps of the inner integrator have to be computed, to find $f\projnkL$ and $f\projnkLp$ from $f\projnk$. These inner steps are, once again, given by the forward Euler method
\begin{align}\label{eq:innerTIntegrator}
	f\projnklp = \FE_{\Dt_0}[f\projnkl] \coloneqq f\projnkl + \Dt_0 \D\e[f\projnkl],
\end{align}
and can be computed directly. Since the outermost extrapolation is once again performed with a first-order approximation of the derivative, this is known as the (two-level) \textit{telescopic projective forward Euler method}.

A schematic summary of the method is presented in \cref{alg:projectiveIntegration}. Once again, the choice of parameters $K_0$, $K_1$, $\Dt_0$, $\Dt_1$, and $\Dt_2$ will be discussed in \cref{sec:spectra}.
\begin{figure}[ht]
	\centering
	\begin{minipage}{.6\linewidth}
		\begin{algorithm}[H]
			\SetAlgoNoLine
			\KwData{$N$, $K_0$, $K_1$, $\Dt_2$, $\Dt_1$, $\Dt_0$}
			$f^0 \leftarrow f(0)$ \;
			\For{$0\leq n\leq N-1$}{
				$f\projnzero \leftarrow f\n$ \;
				\For{$0\leq k\leq K_0$}{
					$f\projnkzero \leftarrow f\projnk$ \;
					\For{$0\leq l\leq K_1$}{
						$f\projnklp \leftarrow \FE_{\Dt_0}[f\projnkl]$ \;
					}
					$f\projnkp \leftarrow \PFE_{\Dt_1}[f\projnk]$ \;
				}
				$f\np \leftarrow \PFE_{\Dt_2}[f\n]$ \;
			}
			\caption{Telescopic Integration (two levels)}
			\label{alg:telescopicIntegration}
		\end{algorithm}
	\end{minipage}
\end{figure}

\paragraph{Efficiency of the Method}

The efficiency of telescopic projective integration can be measured as in the previous section. The $n\th$ level reduces the computational cost of a simulation by a factor of $S_n$, and the overall efficiency factor is $S$:
\begin{align}\label{eq:efficiency}
	S_n
	= \frac{\Dt_{n+1} \prt*{\Dt_n}^{-1}}{K_n+1}
	= \frac{M_n+K_n+1}{K_n+1},\quad
	S = \prod_{n=1}^{N} S_n.
\end{align}
\section{Spectra \& Stability}\label{sec:spectra}

\placedfigure{
	\includegraphics{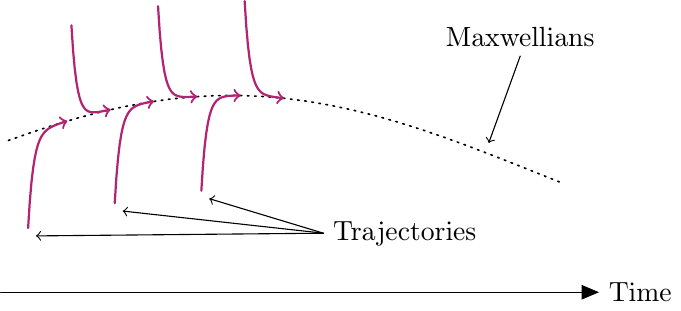}
	\caption{
		Trajectories in solution space for systems with Knudsen number in the hydrodynamic regime.
	}
	\label{fig:manifold}
}
\placedfigure{
	\includegraphics{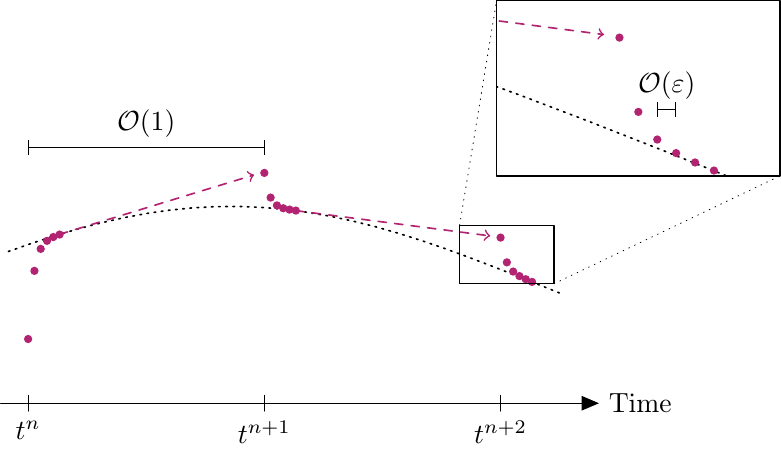}
	\caption{
		Sketch of projective integration: the method performs steps of size $\O(1)$, but introduces intermediate steps of size $\O(\eps)$ to ensure the overall stability of the numerical solution.
	}
	\label{fig:manifoldPoints}
}
\placedfigure{
	\includegraphics{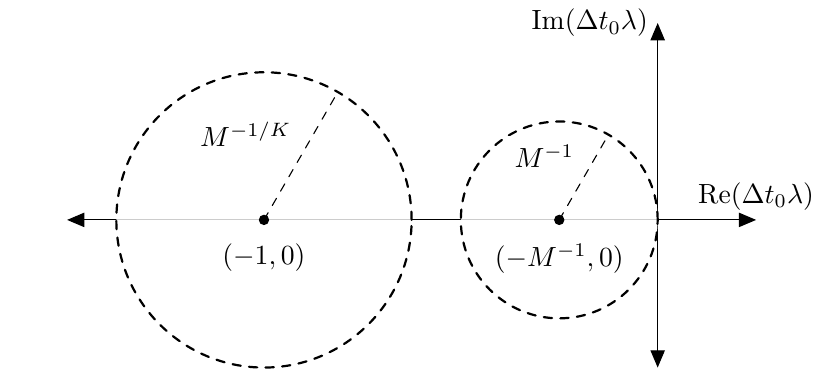}
	\caption{
		Asymptotic stability region for the projective forward Euler integrator.
	}
	\label{fig:PFEStability}
}
\placedfigure{
	\includegraphics{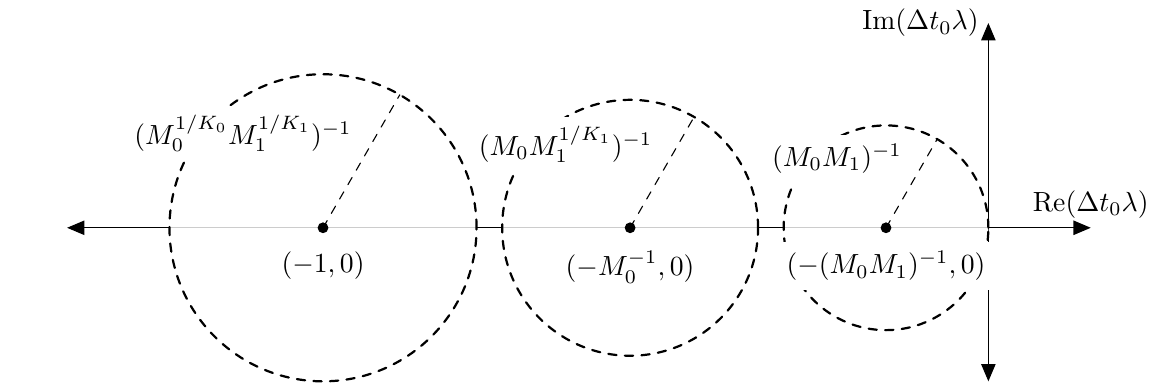}
	\caption{
		Asymptotic stability region for the telescopic projective forward Euler integrator.
	}
	\label{fig:TPFEStability}
}
\placedfigure{
	\includegraphics{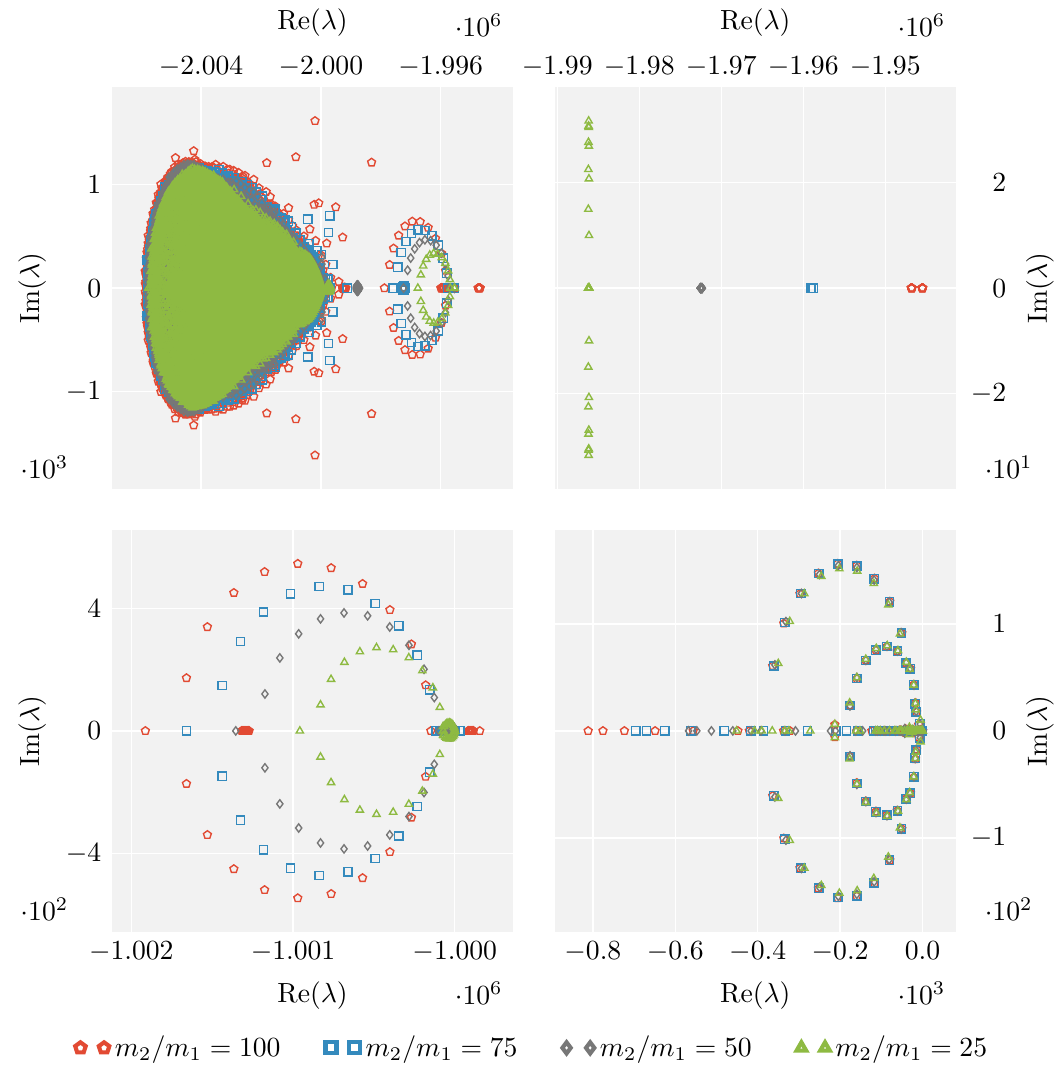}
	\caption{
		Approximate spectra for mixture Sod Tube problem of \cref{sec:ST}: various mass ratios. \textbf{Top}: fast modes, $\O(\eps)$ (split across two regions). \textbf{Bottom left}: middle modes. \textbf{Bottom right}: slow modes, $\O(1)$.
	}
	\label{fig:spectrum}
}

The purpose of this section is to discuss the choice of integration parameters, $\Dt_i$ and $K_i$. We shall recall the stability analysis of \cite{G.K2003,G.K2003a}, and illustrate it in the context of kinetic equations, and the multispecies \ac{BGK} model in particular.

The suitability of projective integration as a numerical method for kinetic equations with very small Knudsen numbers can be understood intuitively by considering their hydrodynamic asymptotics. As referenced in \cref{sec:hydrolimits}, \Cref{eq:kinetic,eq:kineticMultispecies} can be studied in the limit $\varepsilon\rightarrow 0$, and their limiting behaviour can be characterised by systems of equations for the evolution of their moments. In the space of solutions, as the collision operator drives the state towards local equilibrium, trajectories rapidly arrive at a hypersurface comprised only of Maxwellian distributions, see \cref{fig:manifold} for a diagram. Once on this surface, the only possible evolution is within, which corresponds to the evolution of the distribution's moments as governed by the corresponding macroscopic system. For small values of $\varepsilon$, the arrival to and the evolution within the Maxwellian surface occur in very different time scales.

As discussed in \cref{sec:schemes}, resolving both of these time scales with a typical scheme proves very costly, since a typical integrator will require steps of size \revisionTwo{$\Dt = \mathcal{O}\prt{\varepsilon}$} for stability. However, the two-scale nature of the problem can be exploited to escape this restriction: the integrator is permitted large steps (whose size are limited only by a hyperbolic CFL-like condition), as long as they alternate with smaller steps which ensure the return to the equilibrium hypersurface, as sketched in \cref{fig:manifoldPoints}.

This intuitive portrayal can be formalised by performing error stability analysis on the projective integration and telescopic integration methods. An asymptotic stability result is presented in \cite{G.K2003}: for large $M$, the stability region for the projective Euler method approaches two disjoint disks on the $\Dt_0\lambda$ complex plane; $\mathcal{D}\prt*{-1,M^{-1/K}}$ and $\mathcal{D}\prt*{-M^{-1},M^{-1}}$, where by $\mathcal{D}\prt*{c,r}$ we mean the disk centred at the complex number $c$ with radius $r$ (see \cref{fig:PFEStability}). The two-disk stability region is perfectly suited for a two-cluster spectrum, which often lies behind the fast-slow dynamics described above. Once the spectrum is known, the parameters $\Dt_0$, $K$, and $M$ can be chosen so that the scaled eigenmodes, $\Dt_0\lambda$, lie within the two stable disks: the fast modes will be contained in the left disk, and the slow modes, in the right disk. As a curiosity, note that the right disk corresponds to the stability region of a standard forward Euler method with step $\Dt_0 M$; this is indeed the step size of the projection \eqref{eq:outerIntegrator}.

In more complex settings, a two-disk stability region might not suffice. For instance, the \ac{BGK} model with a density-dependent collision frequency, the full Boltzmann equation, and the multispecies \ac{BGK} model which we consider in this work all lead to more complicated spectra. Fortunately, telescopic projective methods also lead to richer stability regions, and can be used for our purposes. In the same vein as \cite{G.K2003}, we find that the two-level telescopic projective Euler method comprises three disks when $M_0$ and $M_1$ are large: $\mathcal{D}\prt*{-1,M_0^{-1/K_0}M_1^{-1/K_1}}$, $\mathcal{D}\prt*{-M_0^{-1},M_0^{-1}M_1^{-1/K_1}}$, and $\mathcal{D}\prt*{-M_0^{-1}M_1^{-1},M_0^{-1}M_1^{-1}}$ (see \cref{fig:TPFEStability}). As before, the parameters $\Dt_0$, $K_0$, $K_1$, $M_0$, and $M_1$ have to be chosen so that the scaled eigenmodes lie within the stability region. This three-disk configuration is featured in the spectra of the numerical experiments discussed in \cref{sec:SB,sec:KH,sec:RM}.

\subsection{Estimating Spectra}

The parameter choice strategy presented above relies entirely on accurate spectral information for the problems in question. However, finding the spectrum of kinetic equations, or even establishing partial information (such as bounds, number of clusters, or spectral gaps) is a difficult problem. For the single species Boltzmann equation, we refer to the classic paper of \cite{E.P1975}, its $\Lone$ extension \cite{Rey2013}, \revisionTwo{or its $\Ltwo$ extension \cite{Gervais2020}}. The multispecies theory is sparse, though we highlight the recent works \cite{B.D2016,D.J.M+2016}.

Our approach in this work will be to estimate the spectra of the problems numerically. We will consider the semi-discrete problem \eqref{eq:semidiscrete}, evaluate the Jacobian of the operator $\D\e$ through finite differencing, and compute the eigenvalues of the resulting matrix as an approximate spectrum for the problem. This approach proves extremely successful, as it permits finding optimal integration parameters for all of the numerical experiments of \cref{sec:experiments}. Furthermore, exploiting the mild dependence of the spectra in the dimension of the problem \cite{E.P1975}, we are able to use the spectra of one-dimensional problems in order to find suitable parameters for problems in higher dimensions.

We present in \cref{fig:spectrum} the spectra used in the mixture Sod Tube problem of \cref{sec:ST} for mass ratios ranging from $m_1/m_2=25$ to $m_1/m_2=100$. The slow modes and middle modes, shown in the bottom right and bottom left plots respectively, are seen not to vary significantly across the different ratios, and therefore have similar stability properties.

However, this is not the case for the fastest modes. The fast cluster is shown across two plot panels in \cref{fig:spectrum}: the top left panel comprises the majority of fast modes, which vary slightly; the top right panel shows the rightmost fast modes, which vary greatly as the mass ratio increases. Across the different ratios, the width of the fast cluster grows by an order of magnitude. This spreading of the fast cluster, and the associated reduction in the gap, is consistent with the results in \cite{B.D2016}.

In practice, the widening of the cluster leads to different stability requirements: in the corresponding numerical test, the parameters used in the $m_1/m_2=1$ case also work for the $m_1/m_2=100$ case, with the exception of $K_1$; the step number of the second level is increased from $6$ to $14$. This logarithmic scaling of the computational cost with the stiffness of the problem has been reported in \cite{M.R.S2017} in the single-species setting.

\section{Numerical Experiments}
\label{sec:experiments}

In this section, we apply our schemes to various physical scenarios in order to showcase the robustness and versatility of the methods. We will consider one example with $\dimx=\dimv=1$, the mixture Sod Tube, and three examples with $\dimx=\dimv=2$: a shock-bubble interaction, a Kelvin-Helmholtz instability, and a Richtmyer-Meshkov instability. The examples demonstrate the ability of the schemes to deal with complex scenarios, as well as their superior efficiency.

\subsection{Mixture Sod Tube \& Extreme Mass Ratios (1D / 1D)}
\label{sec:ST}

\placedfigureHT{
	\includegraphics{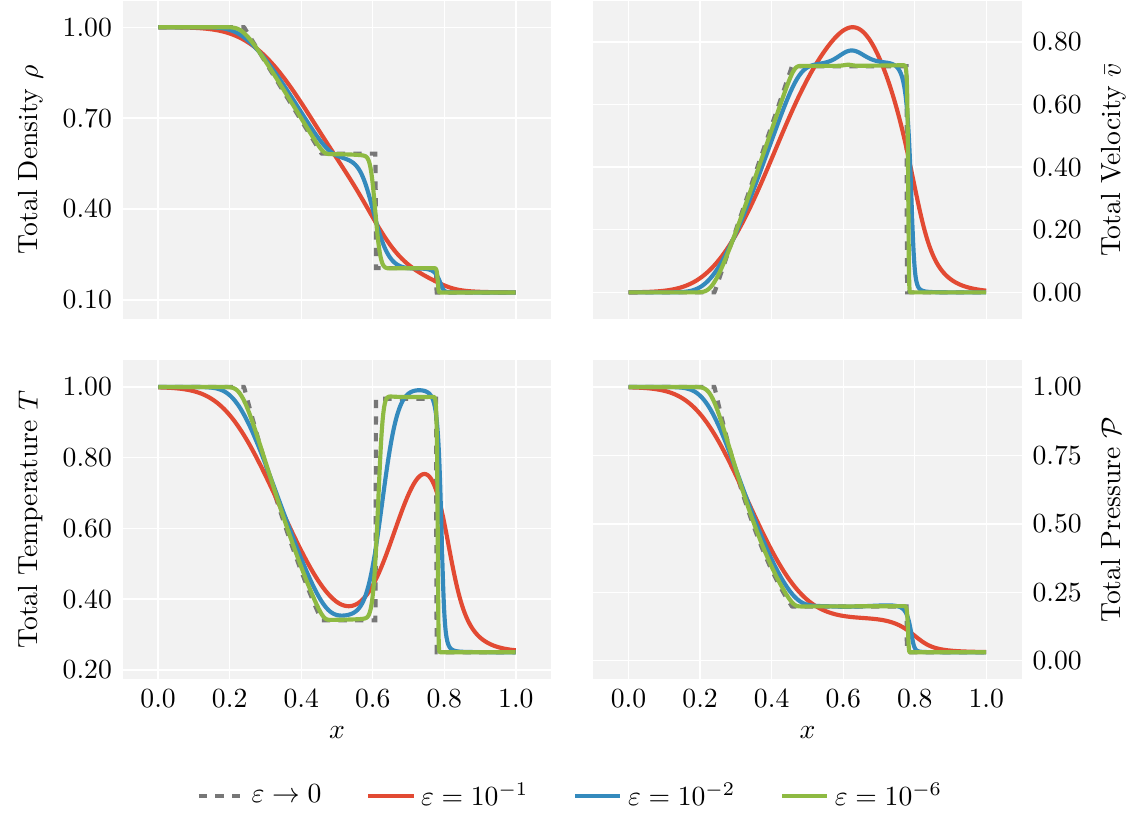}
	\caption{
		Mixture Sod Tube problem: numerical solutions with decreasing Knudsen number and unit mass ratio. The solutions approach the correct $\varepsilon\rightarrow 0$ limit.
	}
	\label{fig:SodTubeEps}
}

\placedfigure{
	\includegraphics{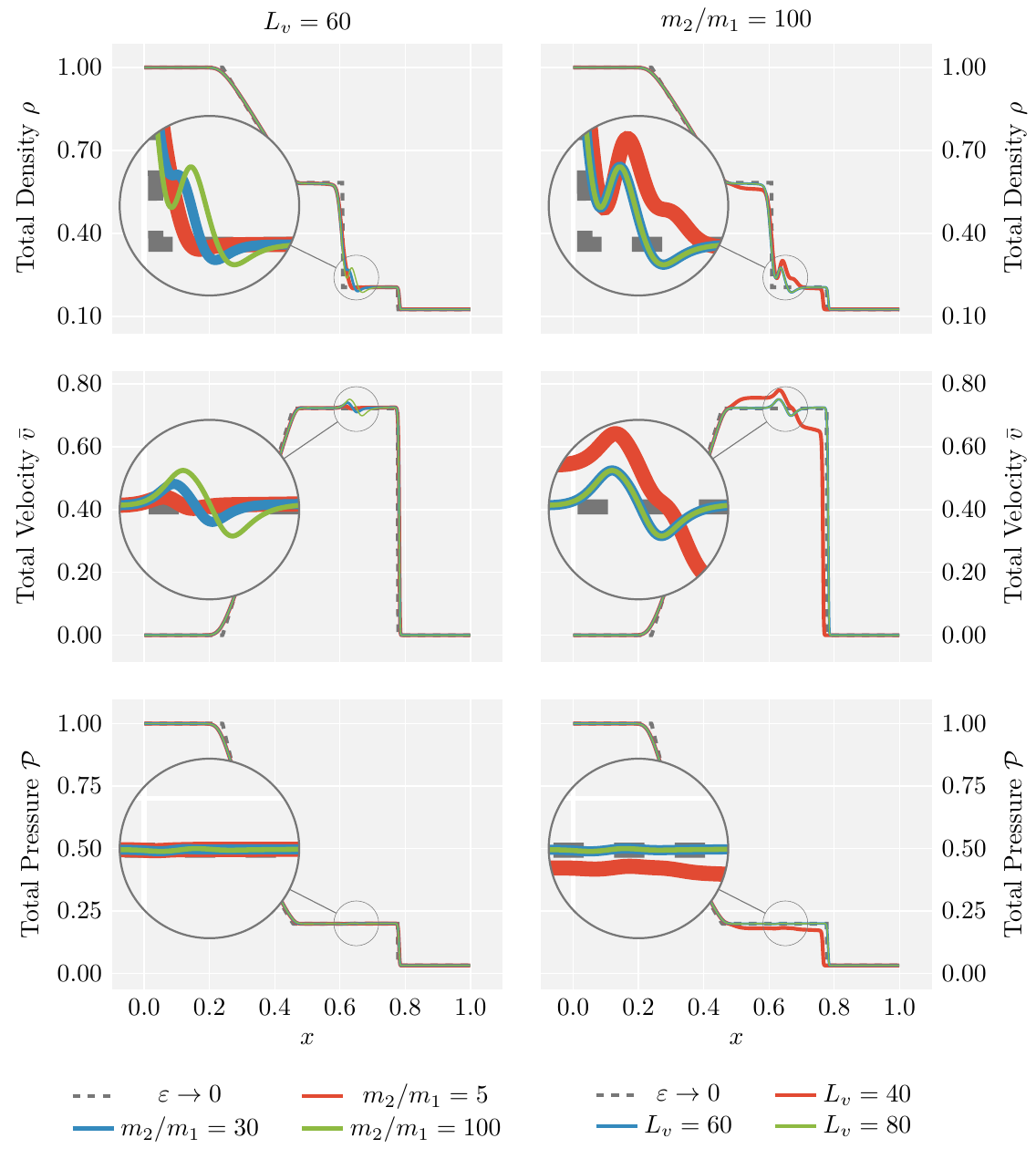}
	\caption{
		Mixture Sod Tube problem: extreme mass ratios. \textbf{Left}: higher mass ratios yield slower hydrodynamic convergence. \textbf{Right}: high mass ratios require larger velocity domains $(-L_v,L_v)$ in order to capture the correct behaviour.
	}
	\label{fig:SodTubeMass}
}

Our first experiment is a generalisation of the Sod Tube problem \cite{Sod1978}. The Sod Tube is a generic Riemann problem for the Euler system \eqref{eq:MultispeciesEuler} (with a single species). A discontinuous initial datum is prescribed, which immediately develops a range of hyperbolic phenomena: a rarefaction wave, a shock wave, and a contact discontinuity. Nevertheless, the solution is analytically tractable, yielding an ideal test case for numerical schemes.

We shall employ here a setting twice removed from the classical Sod Tube. Firstly, we entertain an analogue of the Riemann problem which involves two initially separate gases, whose solution can be constructed in the same fashion as that of the classical problem; see the appendix for the derivation. Secondly, rather than solving a multifluid Euler system, we solve the corresponding problem for the multispecies BGK model, and compare its moments to the aforementioned analytic solution, exploiting the limiting behaviour discussed in \cref{sec:hydrolimits}.

\subsubsection{\revisionOne{Mixture Sod Tube}}

We will consider two gases in a one-dimensional domain $\Omega = (0,1)$. The mass ratio $m_2/m_1$ and Knudsen number will vary through the examples. The initial configuration will emulate a Riemann problem with left state $(\rho_L, v_L, \curlyP_L) = (1,0,1)$ and right state {$(\rho_R, v_R, \curlyP_R) = (2^{-3},0,2^{-5})$}, where the left state is entirely comprised of the first gas, and the right state of the second. To that end, we prescribe their initial distributions as Maxwellians with moments
\begin{equation}\label{eq:SodDatumMaxwellian}
	\left\{
	\begin{array}{lllll}
		\displaystyle
		\rho_1 = (1-\delta) \rho_L,        &
		\rho_2 = \delta \rho_L,            &
		\barv_1 = \barv_2 = v_L,           &
		\curlyP_1 = \curlyP_2 = \curlyP_L, &
		\text{if } x \leq 0.5;
		\\\displaystyle
		\rho_1 = \delta \rho_R,            &
		\rho_2 = (1-\delta) \rho_R,        &
		\barv_1 = \barv_2 = v_R,           &
		\curlyP_1 = \curlyP_2 = \curlyP_R, &
		\text{if } x > 0.5.
	\end{array}
	\right.
\end{equation}
Ideally, we would set $\delta=0$, but this leads to an ill-defined temperature, so we let $\delta=10^{-5}$. This pattern will also be used in later experiments whenever a density is zero. We henceforth refer to this problem as the \textit{mixture Sod Tube problem}.

To begin, we verify the behaviour of this problem in the hydrodynamic limit. We will compute numerically the solution with datum \cref{eq:SodDatumMaxwellian} with mass ratio $m_2/m_1=1$, for three different values of the Knudsen number: $\varepsilon=10^{-1}$, $10^{-2}$, and $10^{-6}$. The larger values of epsilon do not pose significant stiffness, so a direct integration method can be used; for $\varepsilon=10^{-6}$, we resort to a telescopic two-level method.

The solution is computed over the time interval $t\in(0,0.15)$. The domain $\Omega$ is discretised with $\Dx=2^{-10}$, and the velocity space is set as $(-20,20)$, with $\Dv=2^{-4}$. For the direct method employed on the larger values of $\varepsilon$, we let $\Delta t=\num{1.53e-05}$. In the telescopic method, we choose $\Dt_0=\num{5e-07}$, $\Dt_1=\num{2e-06}$, and $\Dt_2=\num{6.1e-05}$, and step numbers $K_0=1$ and $K_1=6$. We impose no-flux boundary conditions.

\Cref{fig:SodTubeEps} shows the numerical solutions superimposed on the analytical limiting solution (computed through the procedure detailed in the appendix). As the Knudsen number decreases, the moments of the solution approach the correct hydrodynamic limit.

\subsubsection{\revisionOne{A CPU Benchmark}}
\label{sec:Benchmark}

\begin{table}\centering

	\begin{tabular}{
			p{0.05\textwidth}
			|
			p{0.07\textwidth}
			p{0.11\textwidth}
			|
			p{0.11\textwidth}
			p{0.12\textwidth}
		}
		\multicolumn{1}{c}{$\varepsilon$}       &
		\multicolumn{2}{c}{CPU Time ($\si{s}$)} &
		\multicolumn{2}{c}{Improvement factor}
		\\
		\midrule
		                                        & Direct
		                                        & Telescopic
		                                        & Real
		                                        & Theoretical
		\\
		$\num{1e-5}$
		                                        & $54$
		                                        & $47$
		                                        & $\num[round-mode=figures,round-precision=3]{1.14893617}$
		                                        & $\num[round-mode=figures,round-precision=3]{3.48772321}$
		\\
		$\num{1e-6}$
		                                        & $508$
		                                        & $47$
		                                        & $\num[round-mode=figures,round-precision=3]{10.808510638}$
		                                        & $\num[round-mode=figures,round-precision=3]{34.877232142857146}$
		\\
		$\num{1e-7}$
		                                        & $5660$
		                                        & $48$
		                                        & $\num[round-mode=figures,round-precision=3]{117.916666667}$
		                                        & $\num[round-mode=figures,round-precision=3]{348.77232142857144}$
		\\
	\end{tabular}

	\caption{
		CPU Time (rounded to nearest second) and improvement factor for the benchmark of \cref{sec:Benchmark}
	}

	\label{tab:CPU}
\end{table}
\revisionOne{
We now revisit the problem \eqref{eq:SodDatumMaxwellian} to conduct a CPU benchmark. We will solve the problem for $\varepsilon=10^{-5}$, $\varepsilon=10^{-6}$, and $\varepsilon=10^{-7}$ with the telescopic method as well as the direct method, and compare the
}

\revisionOne{
The domain $\Omega$ is discretised with $\Dx=2^{-8}$, and the velocity space is set as $(-20,20)$, with $\Dv=2^{-4}$. For the direct method, we let $\Delta t=\varepsilon/2$. In the telescopic method, we choose $\Dt_0=\varepsilon/2$, $\Dt_1=2\varepsilon$, and $\Dt_2=\Dx/16$, and step numbers $K_0=1$ and $K_1=6$. We impose no-flux boundary conditions.
}

\revisionOne{
\cref{tab:CPU} shows the CPU time of each simulation, the real improvement factor (the ratio of CPU times of the direct method to the telescopic method), and the theoretical improvement factor as defined in \cref{eq:efficiency}. The theoretical factor is overly optimistic, as it assumes that the projection steps have negligible computational cost; nevertheless, the telescopic method is faster for $\varepsilon=10^{-5}$ and vastly superior for smaller values of the Knudsen number.
}
\subsubsection{\revisionOne{Extreme Mass Ratios}}

It is of great interest to attempt the numerical solution of problems with extremely large mass ratios. Indeed, a simple mixture of \ch{Ar} and \ch{He} exhibits a mass ratio of $10$, which can be a problem for some numerical methods, as noted in \cite{W.Z.R+2015}. Their method can deal with mass ratios up to $35$; however, ratios twice as large can easily be found in other scenarios, such as mixtures of \ch{H2} and \ch{Xe}.

We will demonstrate the behaviour of our scheme in the context of a hydrodynamic limit under an extreme mass ratio. Conveniently, the asymptotic behaviour of the mixture Sod Tube problem remains unchanged if we assume that the gases have different molecular masses; therefore, it remains a suitable validation case.

We solve the problem for mass ratios $m_2/m_1=5$, $30$, and $100$, in the hydrodynamic regime and compare the solutions. We let $\varepsilon=10^{-6}$. The spatial discretisation is done as before. The velocity space is set as $(-60,60)$, with $\Dv=2^{-4}$; such large domain is unnecessary for the smaller mass ratios, but will be required for $m_2/m_1=100$. We again choose $\Dt_0=\num{5e-07}$, $\Dt_1=\num{2e-06}$, and $\Dt_2=\num{6.1e-05}$, though this time we set $K_0=1$ and $K_1=14$; again, the large step number is only required for the larger mass ratios. In order to justify our velocity discretisation, we will also solve the $m_2/m_1=100$ case with varying velocity spaces $(-L_v,L_v)$, for $L_v=40$, $60$, and $80$, keeping the rest of the parameters fixed.

In all cases we deal with the large velocity supports directly. The recent work \cite{B.P2020} has used a rescaling velocity approach reminiscent of \cite{F.R2013} to overcome the same issue, but their strategy remains limited to mass ratios up to $20$.

\Cref{fig:SodTubeMass} shows the numerical solutions, again superimposed on the limiting solution. The effects of the extreme mass ratios can be seen at point of contact discontinuity (the boundary between the two gases), magnified in the figure. The left column shows the solutions with various mass ratios; the overall hydrodynamic limit is captured well. However, the interfacial effects are more pronounced as the mass ratio increases, and will require a smaller Knudsen number before they become imperceptible. The right column shows the effect of the choice of velocity domain in the solution: $L_v=40$ leads to widespread error, whereas $L_v=60$ recovers the correct behaviour, and is in fact indistinguishable from $L_v=80$.
\subsection{Shock-Bubble Interaction (2D / 2D)}
\label{sec:SB}

\begin{table}\centering

	\begin{tabular}{
			p{0.04\textwidth}
			p{0.05\textwidth}
			|
			p{0.04\textwidth}
			p{0.115\textwidth}
			p{0.08\textwidth}
			|
			p{0.05\textwidth}
			p{0.125\textwidth}
			p{0.115\textwidth}
		}
		\multicolumn{2}{l}{Physical}                   &
		\multicolumn{3}{l}{Phase space discretisation} &
		\multicolumn{3}{l}{Time discretisation}
		\\
		\midrule
		$m_1$                                          & $1$             &
		$\Dx$                                          & $\num{0.025}$   &                    &
		$K_0$                                          & $1$             &
		\\
		$m_2$                                          & $5$             &
		$\Dy$                                          & $\num{0.025}$   & ($=\Dx$)           &
		$K_1$                                          & $6$             &
		\\
		$\varepsilon$                                  & $\num{1e-5}$    &
		$\Dv$                                          & $\num{0.25}$    &                    &
		$\Dt_0$                                        & $\num{0.5e-5}$  & ($=\varepsilon/2$)
		\\
		                                               &                 &
		                                               &                 &                    &
		$\Dt_1$                                        & $\num{2e-5}$    & ($=2\varepsilon$)
		\\
		                                               &                 &
		                                               &                 &                    &
		$\Dt_2$                                        & $\num{0.00125}$ & ($=\Dx/20$)
		\\
	\end{tabular}

	\caption{
		Parameters for the shock-bubble interaction test of \cref{sec:SB}, see \cref{fig:shockBubble2DCmaps,fig:shockBubble2DBubble,fig:shockBubble2DContour}
	}

	\label{tab:shockBubble2DPars}
\end{table}
\placedfigure{
	\includegraphics{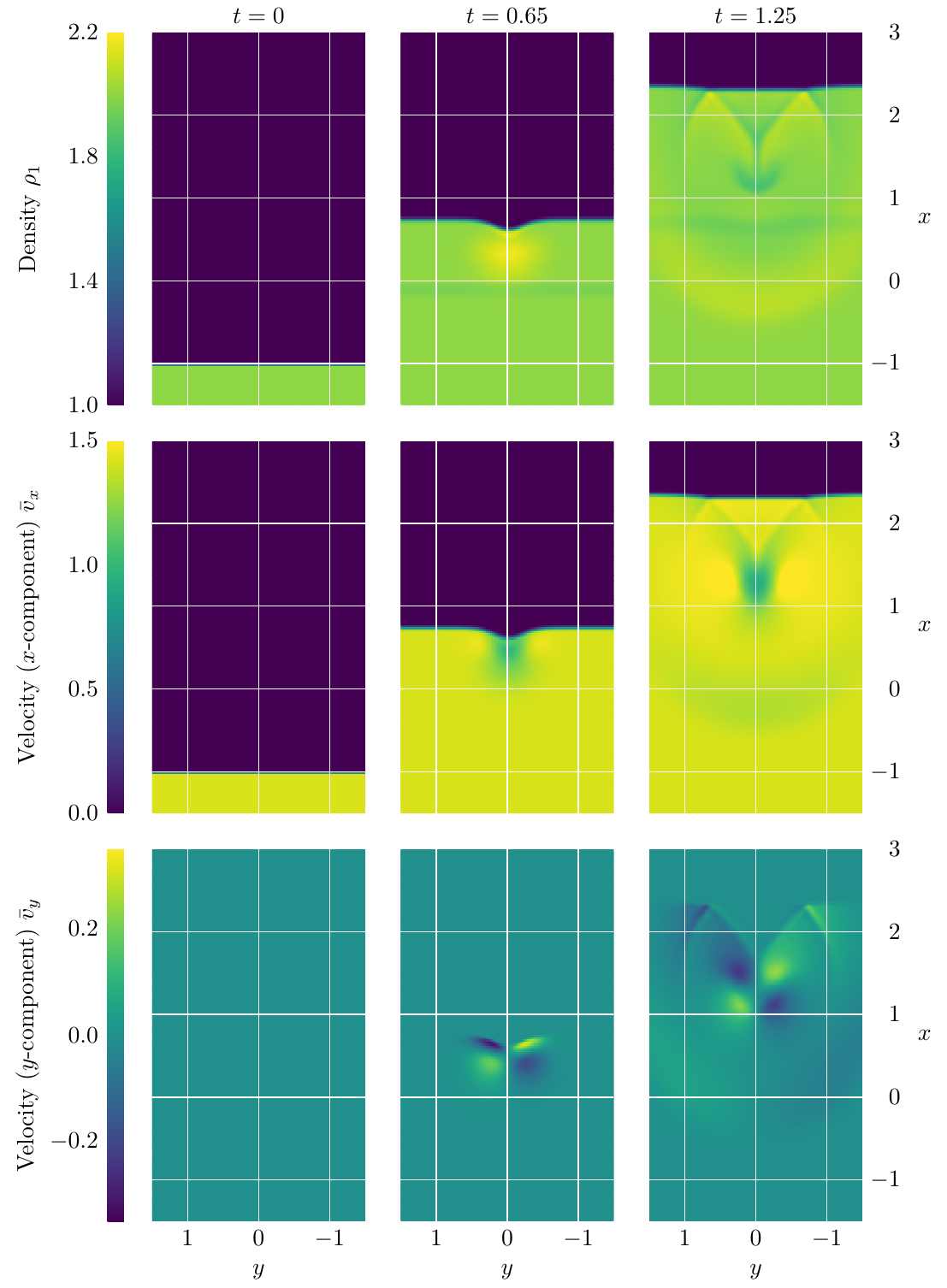}
	\caption{
		Shock-bubble interaction: evolution of the first (lighter) species.
		\textbf{Left}: datum.
		\textbf{Middle}: shock interacts with bubble.
		\textbf{Right}: pressure waves arise as a result of the shock-bubble interaction.
		Simulation parameters can be found in \cref{tab:shockBubble2DPars}.
	}
	\label{fig:shockBubble2DCmaps}
}

\placedfigure{
	\includegraphics{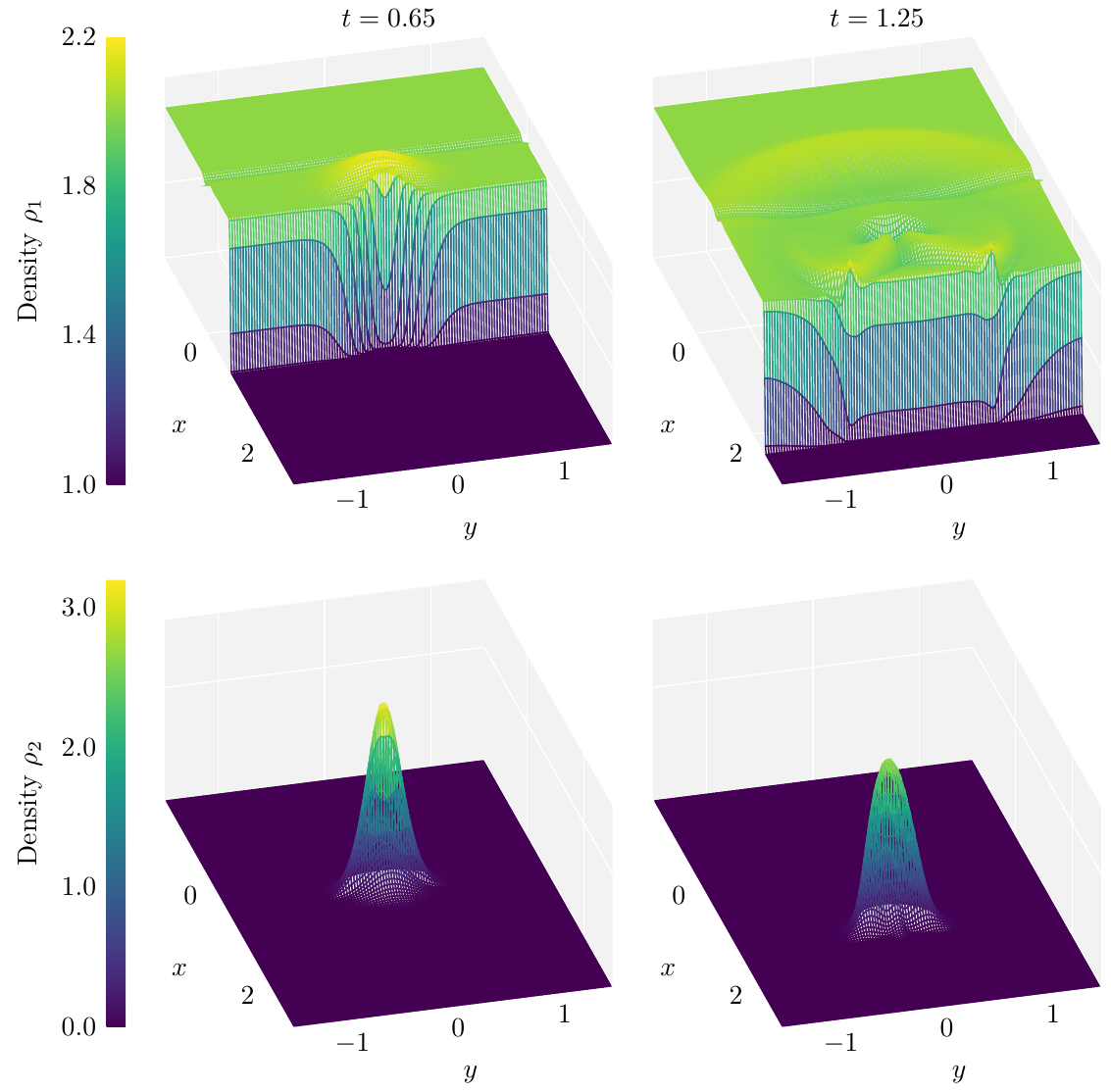}
	\caption{
		Shock-bubble interaction: evolution of the densities.
		\textbf{Left}: shock interacts with bubble.
		\textbf{Right}: bubble is dragged and deformed by the stream.
		Simulation parameters can be found in \cref{tab:shockBubble2DPars}.
	}
	\label{fig:shockBubble2DBubble}
}

\placedfigure{
	\includegraphics{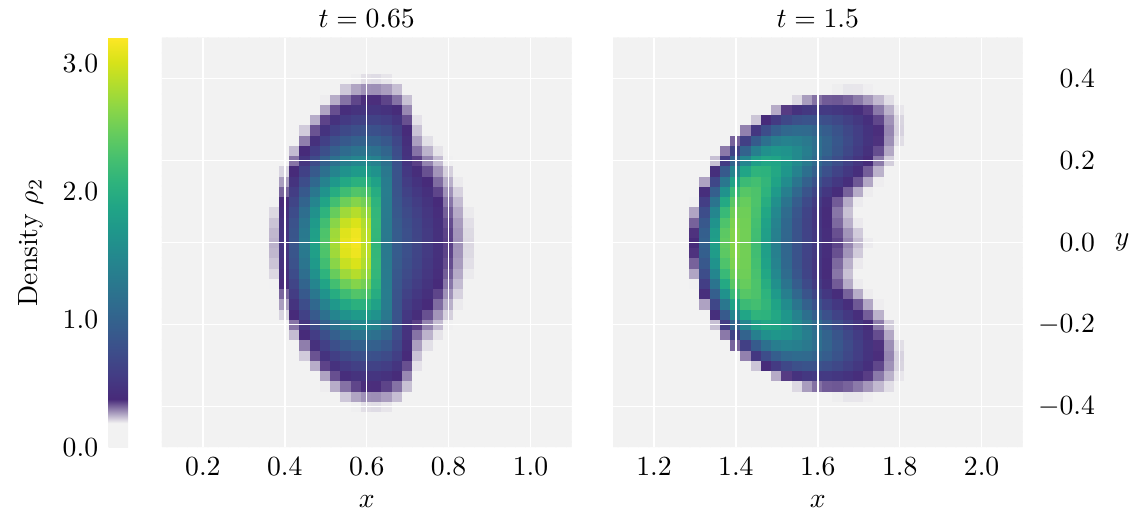}
	\caption{
		Shock-bubble interaction: detail of the bubble (second species).
		\textbf{Left}: the bubble has been compressed by the shock.
		\textbf{Right}: the bubble is displaced by the stream and deforms into a ``C shape''.
		Simulation parameters can be found in \cref{tab:shockBubble2DPars}.
	}
	\label{fig:shockBubble2DContour}
}

We investigate the interaction between a travelling shock and a smooth stationary bubble. This is a multispecies adaptation of a one-species test case proposed in \cite{Torrilhon2006}; the original test has been used to validate numerical schemes for the \ac{BGK} model \cite{C.L2010}, including projective integration \cite{M.R.S2019}.

We will consider a mixture of two gases with mass ratio $m_2/m_1=5$ in a rectangular domain $\Omega = (-1.5,3)\times(-1.5,1.5)$. The Knudsen number is chosen in the hydrodynamic regime, $\varepsilon=\num{1e-5}$. The initial configuration of the first (lighter) gas is a normal shock wave of Mach number $\Mach=2$ propagating in the positive $x$-direction. Its initial distribution, $f_1(0,\bx,\bv)$, is chosen as the Maxwellian corresponding to the following Riemann datum:
\begin{equation}\label{eq:SBDatum1}
	\left\{
	\begin{array}{lllll}
		\displaystyle
		\rho_1 = 2,          &
		\barv_{1,x} = 1.414, &
		\barv_{1,y} = 0,     &
		T_1 = 2.5,           &
		\text{if } x \leq -1;
		\\\displaystyle
		\rho_1 = 1,          &
		\barv_{1,x} = 0,     &
		\barv_{1,y} = 0,     &
		T_1 = 1,             &
		\text{if } x > -1.
	\end{array}
	\right.
\end{equation}
The second (heavier) gas is initially stationary: a Maxwellian with density
\begin{align}\label{eq:SBDatum2}
	\rho_2(0,\bx) = \exp\set*{-16\abs*{\bx-\bx_0}^2}.
\end{align}
The temperatures of both gases around the bubble are chosen equal to each other and in such a way that there is unit total pressure:
\begin{align}\label{eq:SBDatum3}
	T_1(0,\bx) = T_2(0,\bx) = \prt*{n_1(0,\bx) + n_2(0,\bx)}^{-1},
\end{align}
which is consistent with the right state in \cref{eq:SBDatum1}.

Numerically, the domain $\Omega$ is discretised with $\Dx=\Dy=\num{0.025}$ ($180\times 120$ cells). The velocity space is set as $(-12,12)^2$, with $\Dv=\num{0.25}$ ($96^2$ cells). In total, the phase space is discretised with $\num[group-separator = {,},scientific-notation=false]{199065600}$ cells. We impose outflow boundary conditions.

The solution corresponding to \cref{eq:SBDatum1,eq:SBDatum2,eq:SBDatum3} is computed over the time $t\in(0,1.5)$. We employ the telescopic two-level method, with step sizes $\Dt_0=\num{0.5e-5}$, $\Dt_1=\num{2e-5}$, and $\Dt_2=\num{0.00125}$, and step numbers $K_0=1$ and $K_1=6$. The inner steps follow the pattern $4\Dt_0 = 2\varepsilon = \Dt_1$. The outermost step is restricted by the stability of the transport scheme, rather than the projective integration. The use of telescopic projective integration decreases the computational cost by a factor of $S\simeq 18$, as defined in \cref{eq:efficiency}. This factor, considering the fine and high-dimensional mesh employed here, is extremely beneficial, reducing the simulation time from a matter of days to a matter of hours!

\Cref{fig:shockBubble2DCmaps} shows the evolution of the first (lighter) gas at several times. The shock travels to meet the bubble; when they meet, the central ($y\simeq 0$) portion of the shock is slowed down. As the shock traverses the bubble, two pressure waves are formed: one travelling upstream, appears as a reflection from the initial shock-bubble interaction; the second, travelling downstream, arises as the shock surrounds the obstacle and both ``arms'' meet behind the bubble.

\Cref{fig:shockBubble2DBubble} compares the densities of either gas, and \cref{fig:shockBubble2DContour} presents a detailed view of the bubble. The bubble is compressed by the shock when it first arrives. As the shock passes, the bubble is displaced downstream and deforms into a ``C shape''.

A summary of the simulation parameters can be found in \cref{tab:shockBubble2DPars}.

\subsection{Kelvin-Helmholtz Instability (2D / 2D)}
\label{sec:KH}

\begin{table}\centering

	\begin{tabular}{
			p{0.04\textwidth}
			p{0.05\textwidth}
			|
			p{0.04\textwidth}
			p{0.145\textwidth}
			p{0.085\textwidth}
			|
			p{0.05\textwidth}
			p{0.160\textwidth}
			p{0.115\textwidth}
		}
		\multicolumn{2}{l}{Physical}                   &
		\multicolumn{3}{l}{Phase space discretisation} &
		\multicolumn{3}{l}{Time discretisation}
		\\
		\midrule
		$m_1$                                          & $1$                 &
		$\Dx$                                          & $\num{0.0078125}$   & ($=2^{-7}$)        &
		$K_0$                                          & $2$                 &
		\\
		$m_2$                                          & $5$                 &
		$\Dy$                                          & $\num{0.0078125}$   & ($=\Dx$)           &
		$K_1$                                          & $4$                 &
		\\
		$\varepsilon$                                  & $\num{1e-5}$        &
		$\Dv$                                          & $\num{0.5}$         &                    &
		$\Dt_0$                                        & $\num{0.5e-5}$      & ($=\varepsilon/2$)
		\\
		                                               &                     &
		                                               &                     &                    &
		$\Dt_1$                                        & $\num{2e-5}$        & ($=2\varepsilon$)
		\\
		                                               &                     &
		                                               &                     &                    &
		$\Dt_2$                                        & $\num{0.000390625}$ & ($=\Dx/20$)
		\\
	\end{tabular}

	\caption{
		Parameters for the Kelvin-Helmholtz instability test of \cref{sec:KH}.
	}

	\label{tab:KelvinHelmholtz2DPars}
\end{table}

\placedfigure{
	\includegraphics{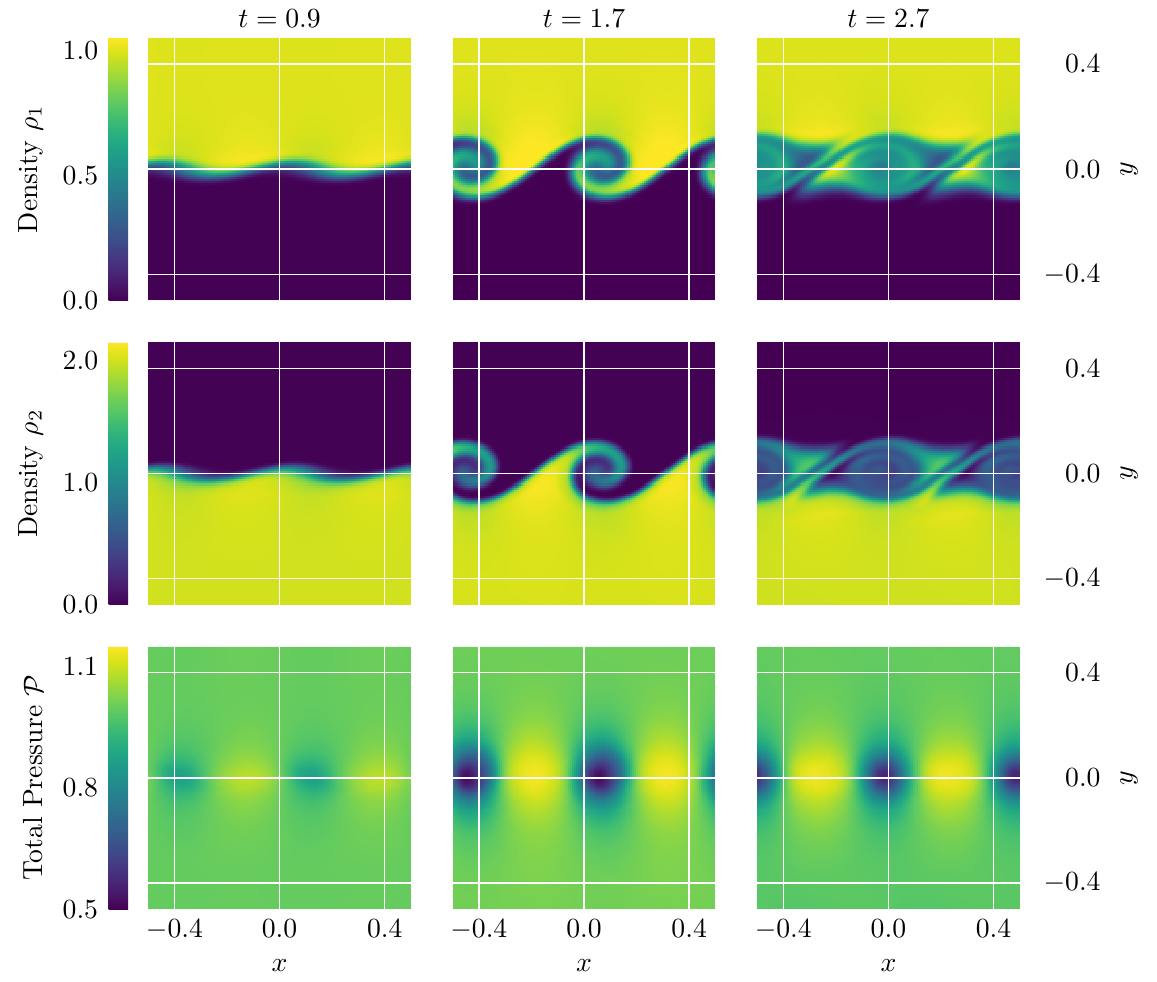}
	\caption{
		Kelvin-Helmholtz instability: evolution in time.
		\textbf{Left}: undulations appear on the gas interface.
		\textbf{Middle}: clear vortices have developed.
		\textbf{Right}: vortices begin to smear and merge.
		Simulation parameters can be found in \cref{tab:KelvinHelmholtz2DPars}.
	}
	\label{fig:KelvinHelmholtz2DCmapsThree}
}

\placedfigure{
	\includegraphics{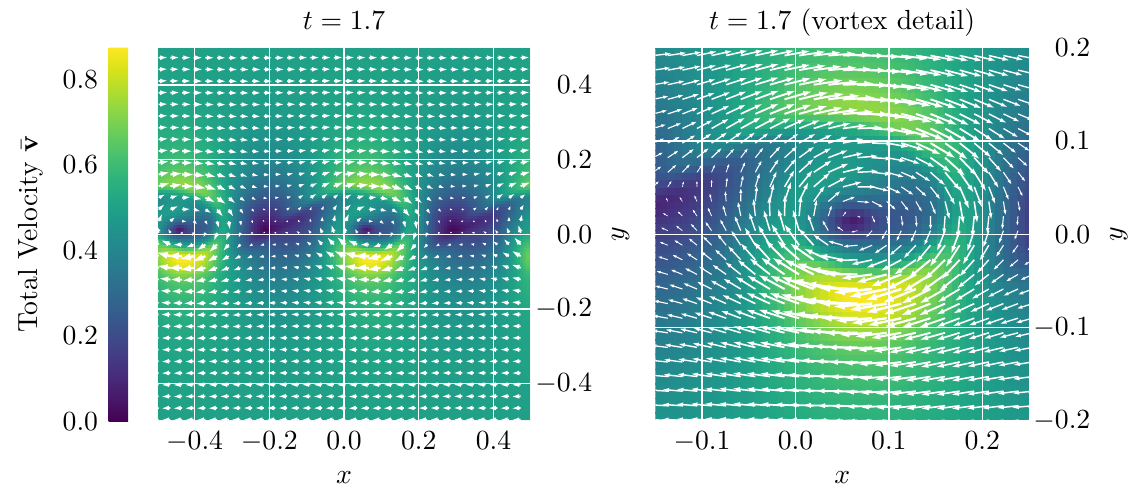}
	\caption{
		Kelvin-Helmholtz instability: modulus of the total velocity $\barbv$ and stream lines at peak vorticity.
		\textbf{Left}: whole domain.
		\textbf{Right}: detail of a vortex.
		Simulation parameters can be found in \cref{tab:KelvinHelmholtz2DPars}.
	}
	\label{fig:KelvinHelmholtz2DStream}
}

We turn our attention to the Kelvin-Helmholtz instability. This is a well-known phenomenon where vortices arise at the interface of two fluids with different density moving at different speeds. This test was also used to validate projective integration for the single-species \ac{BGK} model in \cite{M.R.S2019}, with an initial configuration drawn from \cite{M.L.P2012}. Here we shall present a two-species analogue.

We consider two gases with mass ratio $m_2/m_1=5$ in a square domain $\Omega = (-0.5,0.5)^2$, periodic along the $x$-direction. The Knudsen number is chosen in the hydrodynamic regime, $\varepsilon=\num{1e-5}$; higher values do not yield sufficiently defined vortices. The gases are initially separate: the first (lighter) gas occupies the $y\geq 0$ region, with density $\rho_1=1$ and horizontal velocity $\barv_{1,x}=0.5$; the second (heavier) gas occupies the $y< 0$ region, with density $\rho_2=2$ and horizontal velocity $\barv_{2,x}=-0.5$. Their temperatures are chosen to ensure there is equal unit pressure across the boundary: $T_1=n_1^{-1}$, $T_2=n_2^{-1}$. Both gases are given a small vertical velocity to induce vorticity:
$
	\barv_{1,y} = \barv_{2,y} = \num{1e-2} \sin(4\pi x).
$
The initial distribution of the gases are the Maxwellians corresponding to these moments.

Numerically, the domain $\Omega$ is discretised with $\Dx=\Dy=\num{0.0078125}$ ($128^2$ cells). The velocity space is set as $(-8,8)^2$, with $\Dv=\num{0.5}$ ($32^2$ cells). In total, the phase space is discretised with $\num[group-separator = {,},scientific-notation=false]{16777216}$ cells. We impose outflow boundary conditions along the non-periodic boundaries.

The solution is computed over the time $t\in(0,3)$. We employ the telescopic two-level method, with step sizes $\Dt_0=\num{0.5e-5}$, $\Dt_1=\num{2e-5}$, and $\Dt_2=\num{0.000390625}$, and inner steps $K_0=2$ and $K_1=4$. Again the inner steps follow the pattern $4\Dt_0 = 2\varepsilon = \Dt_1$, and the outermost step is restricted by the hyperbolic CFL condition. The use of telescopic projective integration decreases the computational cost by a factor of $S\simeq 5$, as defined in \cref{eq:efficiency}.

\Cref{fig:KelvinHelmholtz2DCmapsThree} shows the evolution of the densities of each gas as well as the total pressure. The small initial vertical perturbation results in the undulation of the interface. The relative difference in horizontal velocity rapidly causes the formation of clearly defined vortices. As the shear stress, the vortices are smeared horizontally.

\Cref{fig:KelvinHelmholtz2DStream} studies the vortex structure at the time corresponding to the middle column of \cref{fig:KelvinHelmholtz2DCmapsThree}, where vorticity is at a peak. The figure displays the modulus of the total velocity $\barbv$; the streamlines are superimposed.

A summary of the simulation parameters can be found in \cref{tab:KelvinHelmholtz2DPars}.

\subsection{Richtmyer-Meshkov Instability (2D / 2D)}
\label{sec:RM}

\begin{table}\centering

	\begin{tabular}{
			p{0.04\textwidth}
			p{0.05\textwidth}
			|
			p{0.04\textwidth}
			p{0.110\textwidth}
			p{0.085\textwidth}
			|
			p{0.05\textwidth}
			p{0.125\textwidth}
			p{0.115\textwidth}
		}
		\multicolumn{2}{l}{Physical}                   &
		\multicolumn{3}{l}{Phase space discretisation} &
		\multicolumn{3}{l}{Time discretisation}
		\\
		\midrule
		$m_1$                                          & $1$             &
		$\Dx$                                          & $\num{0.00025}$ &                    &
		$K_0$                                          & $1$             &
		\\
		$m_2$                                          & $5$             &
		$\Dy$                                          & $\num{0.00025}$ & ($=\Dx$)           &
		$K_1$                                          & $6$             &
		\\
		$\varepsilon$                                  & $\num{1e-6}$    &
		$\Dv$                                          & $\num{0.25}$    &                    &
		$\Dt_0$                                        & $\num{0.5e-6}$  & ($=\varepsilon/2$)
		\\
		                                               &                 &
		                                               &                 &                    &
		$\Dt_1$                                        & $\num{2e-6}$    & ($=2\varepsilon$)
		\\
		                                               &                 &
		                                               &                 &                    &
		$\Dt_2$                                        & $\num{6.25e-5}$ & ($=\Dx/4$)
		\\
	\end{tabular}

	\caption{
		Parameters for the Richtmyer-Meshkov instability test of \cref{sec:RM}.
	}

	\label{tab:RichtmyerMeshkov2DPars}
\end{table}

\placedfigure{
	\includegraphics{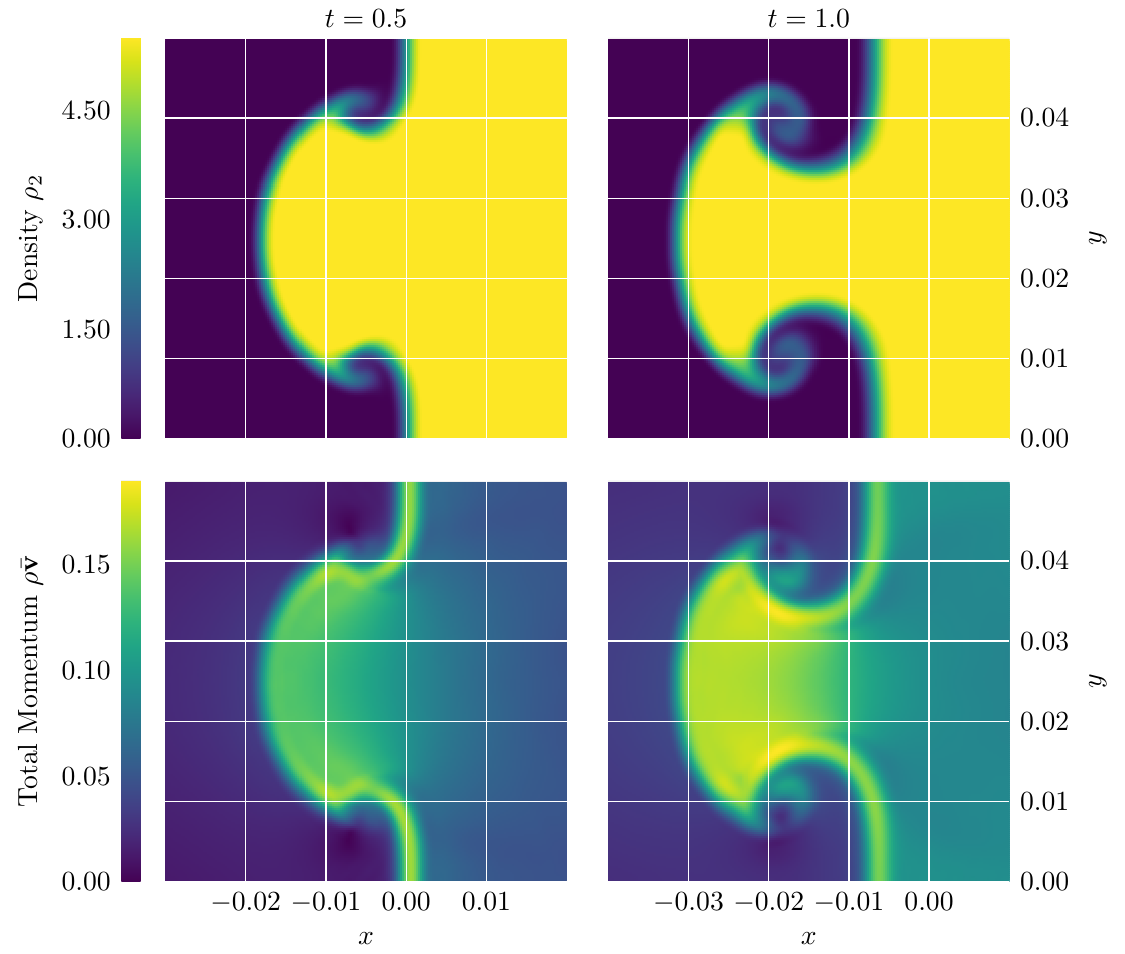}
	\caption{
		Richtmyer-Meshkov instability: density of the second species and modulus of the total momentum.
		\textbf{Left}: after the initial shock interaction, vortices arise.
		\textbf{Right}: the interface is deformed into a mushroom-like shape.
		Simulation parameters can be found in \cref{tab:RichtmyerMeshkov2DPars}.
	}
	\label{fig:RichtmyerMeshkov2DEnd}
}

\placedfigure{
	\includegraphics{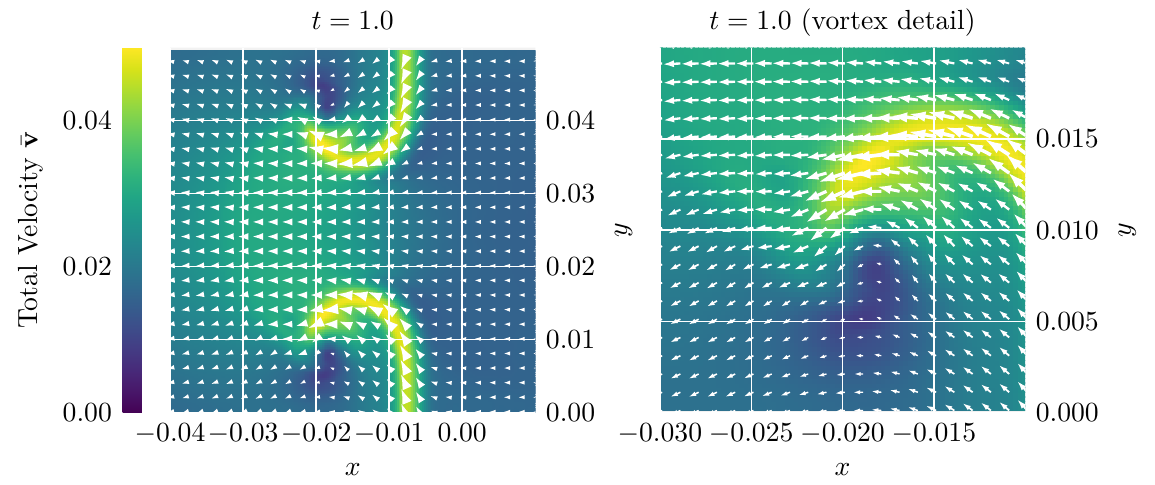}
	\caption{
		Richtmyer-Meshkov instability: modulus of the total velocity $\barbv$ and stream lines at the final time.
		\textbf{Left}: whole domain.
		\textbf{Right}: detail of the lower vortex.
		Simulation parameters can be found in \cref{tab:RichtmyerMeshkov2DPars}.
	}
	\label{fig:RichtmyerMeshkov2DStream}
}

\placedfigure{
	\includegraphics{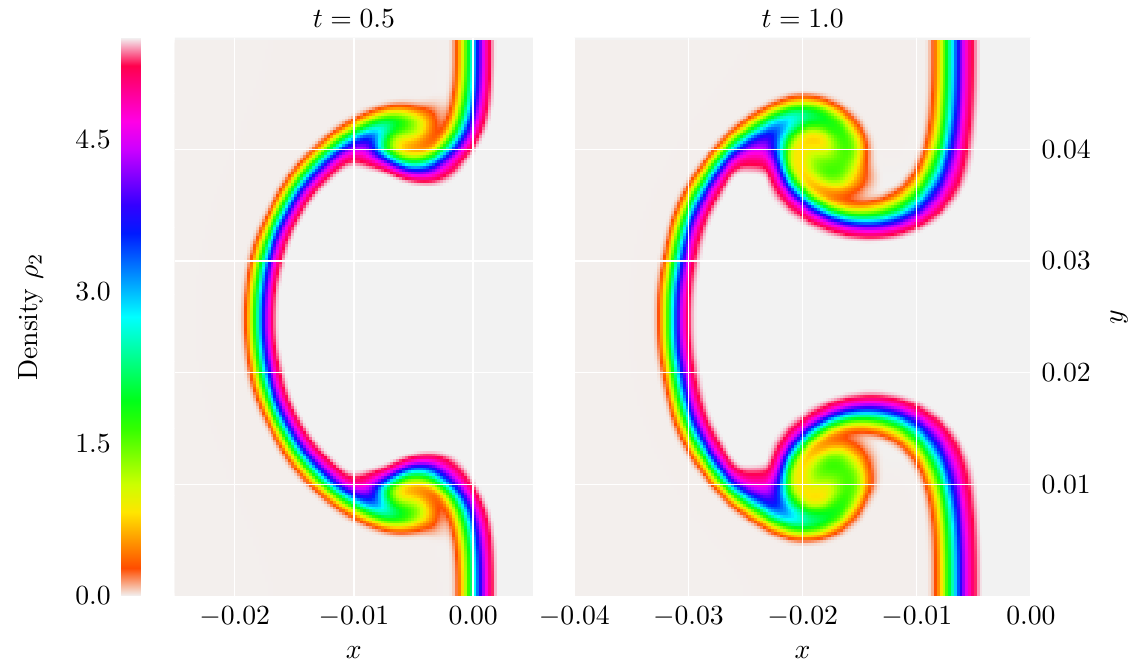}
	\caption{
		Richtmyer-Meshkov instability: detail of the gas interface (density of the second species).
		\textbf{Left}: vortices appear.
		\textbf{Right}: mushroom-like interface is formed.
		Simulation parameters can be found in \cref{tab:RichtmyerMeshkov2DPars}.
	}
	\label{fig:RichtmyerMeshkov2DSequenceShort}
}

To conclude this section, we study the Richtmyer-Meshkov instability. This is a phenomenon which takes place when the perturbed interface between a thin and a dense gas is momentarily accelerated by a passing shock. The misalignment of the pressure gradient (due to the shock) and the density gradient induces vorticity, which leads to the formation of a mushroom-like interface\footnote{This instability is not to be confused with the related Rayleigh-Taylor instability, which arises on a perturbed interface when the dense gas is accelerated continuously by the thin gas, rather than by a shock. More details as well as experimental images can be found in \cite{W.N.J2001}.}. A more detailed physical background, as well as experimental images, can be found in \cite{A.J2008}. This instability has been used to validate a lattice-Boltzmann method in \cite{S.M2020}.

We shall consider two gases with mass ratio $m_2/m_1=5$ in a rectangular domain $\Omega = (-0.5,0.5)\times(0,0.5)$. The Knudsen number is here chosen as $\varepsilon=\num{1e-6}$ in order to correctly resolve the instability. The gases are initially separate: the first (lighter) gas occupies the $x\leq b(y)$ region, with a perturbed boundary given by
$
	x = - \num{1e-2} \sin(20\pi y),
$
whereas the second (heavier) gas lies in the rest of the domain. The initial configuration of the first gas is a normal shock wave of Mach number $\Mach=1.21$ (inspired by \cite{S.M2020}) propagating in the positive $x$-direction. Its initial distribution is chosen according to the following Riemann datum:
\begin{equation}\label{eq:RMDatum1}
	\left\{
	\begin{array}{lllll}
		\displaystyle
		\rho_1 = 1.268,      &
		\barv_{1,x} = 0.256, &
		\barv_{1,y} = 0,     &
		p_1 = 0.809,         &
		\text{if } x \leq s_0;
		\\\displaystyle
		\rho_1 = 1,          &
		\barv_{1,x} = 0,     &
		\barv_{1,y} = 0,     &
		p_1 = 0.5,           &
		\text{if } x > s_0;
	\end{array}
	\right.
\end{equation}
$s_0$, the initial position of the shock, is a negative constant. The second gas is initially stationary, with density $\rho_2=5$, and pressure $p_2=0.5$ to ensure the boundary is initially not forced. The initial distribution of the gases are the Maxwellians corresponding to these moments. For convenience, in order to avoid the interface leaving the domain, all the horizontal velocities are decreased by $\num{0.07}$.

Numerically, the domain $\Omega$ is discretised with $\Dx=\Dy=\num{0.00025}$ ($400\times200$ cells). The velocity space is set as $(-4,4)^2$, with $\Dv=\num{0.25}$ ($32^2$ cells). In total, the phase space is discretised with $\num[group-separator = {,},scientific-notation=false]{81920000}$ cells. We impose outflow boundary conditions.

The solution is computed over the time $t\in(-0.02,1.0)$. We choose the negative initial time and let $s_0=\num{0.0242}$ so that the shock crosses the $x=0$ line exactly at $t=0$. We employ the telescopic two-level method, with step sizes $\Dt_0=\num{0.5e-6}$, $\Dt_1=\num{2e-6}$, and $\Dt_2=\num{6.25e-5}$, and inner steps $K_0=1$ and $K_1=6$. Once more, the inner steps follow the pattern $4\Dt_0 = 2\varepsilon = \Dt_1$ and the outermost step is restricted by the hyperbolic CFL condition. The use of telescopic projective integration decreases the computational cost by a factor of $S\simeq 9$, as defined in \cref{eq:efficiency}.

\Cref{fig:RichtmyerMeshkov2DEnd} shows the density of the heavier gas and the modulus of the total momentum of the system at two different times. The gas interface is seen deforming as it develops a mushroom-like shape. Vortices are visible on either side of the perturbation, and much of the system's momentum is found at the interface.

\Cref{fig:RichtmyerMeshkov2DStream} studies the vortices at the final simulation time. The figure shows the modulus of the total velocity with streamlines, detailing the high speed along the lateral gas interface, and offering a magnified view of the lower vortex.

\Cref{fig:RichtmyerMeshkov2DSequenceShort} offers a detailed view of the gas interface. The density of the heavier gas has been recoloured to omit the extrema, highlighting only the density gradient. The vortex formation and the deformation of the interface are particularly clear here.

A summary of the simulation parameters can be found in \cref{tab:RichtmyerMeshkov2DPars}.

\paragraph*{Acknowledgments}
\addcontentsline{toc}{section}{Acknowledgments}
TR was partially funded by Labex CEMPI (ANR-11-LABX-0007-01) and ANR Project MoHyCon (ANR-17-CE40-0027-01). RB was partially funded by Labex CEMPI (ANR-11-LABX-0007-01) and by the Advanced Grant Nonlocal-CPD (Nonlocal PDEs for Complex Particle Dynamics: Phase Transitions, Patterns and Synchronization) of the European Research Council Executive Agency (ERC) under the European Union's Horizon 2020 research and innovation programme (grant agreement No. 883363).
\singleappendixtitle{Solving the \revisionOne{Mixture Sod Tube Problem}}

\placedfigure{
	\includegraphics{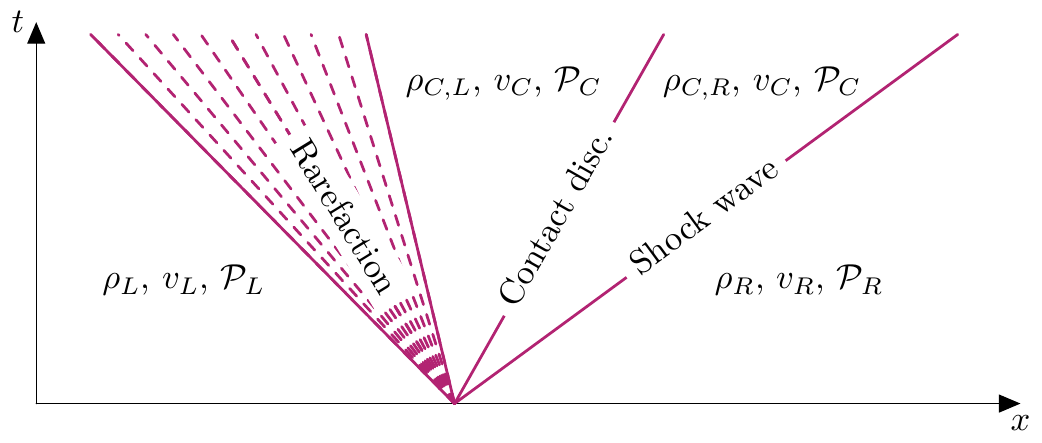}
	\caption{
		Structure of the solution to the Sod Tube problem
	}
	\label{fig:RiemannProblem}
}

We recall, for the sake of a complete exposition, the solution to the classical Sod Tube problem \revisionOne{and the mixture Sod Tube problem}, as can be found in the \revisionOne{literature of hyperbolic problems \cite{Toro1999}}. In the usual setting, the initial Riemann datum is given by the left ($\rho_L$, $v_L$, $\curlyP_L$) and right ($\rho_R$, $v_R$, $\curlyP_R$) states:
\begin{equation}
	\left\{
	\begin{array}{llll}
		\displaystyle
		\rho = \rho_L,       &
		v = v_L,             &
		\curlyP = \curlyP_L, &
		\text{if } x \leq 0;
		\\\displaystyle
		\rho = \rho_R,       &
		v = v_R,             &
		\curlyP = \curlyP_R, &
		\text{if } x > 0.
	\end{array}
	\right.
\end{equation}
As the datum evolves, typically we observe two non-linear waves enveloping a new central region; see \cref{fig:RiemannProblem} for a diagram. This region comprises two states, separated by a contact discontinuity. Both have equal speed and pressure, $v_C$, $\curlyP_C$, but different densities, $\rho_{C,L}$, and $\rho_{C,R}$. We will assume here $\curlyP_L \geq \curlyP_C \geq \curlyP_R$, meaning the leftmost wave is a rarefaction, and the rightmost wave is a shock.

The solution strategy involves first finding values $v_C$, $\curlyP_C$ which are consistent with the rarefaction/shock wave structure on either side. Given a guess for $\curlyP_C$, the value of $v_C$ can be computed from the left state using the \textit{Riemann invariants} (quantities which are constant across the rarefaction). The precise relation is
\begin{align}
	v_{\text{rare}}(\curlyP_C)
	=
	v_L +
	2c_L
	\brk*{1-\prt*{\frac{\curlyP_C}{\curlyP_L}}^{\frac{\gamma-1}{2\gamma}}}
	\prt*{\gamma-1}^{-1},
\end{align}
where $c_L$ is the speed of sound on the left state, computed via $c=\sqrt{\gamma \curlyP/\rho}$, and where $\gamma$ is the \textit{adiabatic exponent}; for an ideal monatomic gas, it is given by $\gamma=1+2/D$, where $D$ is the dimension. Similarly, the value of $v_C$ can be computed from the right state using the \textit{Rankine-Hugoniot conditions}, which relate the quantities on either side of the shock. In this case, the expression is
\begin{align}
	v_{\text{shock}}(\curlyP_C)
	=
	v_R +
	2c_R
	\prt*{1 - \frac{\curlyP_C}{\curlyP_R}}
	\brk*{
		2\gamma \prt*{\gamma-1 + \frac{\curlyP_C\prt{\gamma+1}}{\curlyP_R}}
	}^{-\frac{1}{2}}.
\end{align}
A suitable value of $\curlyP_C$ can thus be found by solving $ v_{\text{rare}}(\curlyP_C) = v_{\text{shock}}(\curlyP_C) $. This determines the value of $v_C$ also. Having established these, the densities on either side of the contact discontinuity are found as
\begin{align}
	\rho_{C,L} = \rho_L
	\prt*{\frac{\curlyP_C}{\curlyP_L}}^{\frac{1}{\gamma}}
	, \quad
	\rho_{C,R} = \rho_R
	\prt*{\frac{\curlyP_C}{\curlyP_R} + \frac{\gamma-1}{\gamma+1}}
	\prt*{\frac{\curlyP_C\prt{\gamma-1}}{\curlyP_R\prt{\gamma+1}} + 1}^{-1}
	.
\end{align}
A detailed derivation of these relations can be found in \cite{Toro1999}, for instance.

We now proceed to establish the spatial structure of the solution. The rarefaction, contact discontinuity, and shock wave all radiate from the initial point of discontinuity. The head and the tail of the rarefaction radiate at speeds $v_L-c_L$ and $v_C-c_C$ respectively. The contact discontinuity travels at speed $v_C$. Lastly, the shock propagates at speed
$
	s = \frac{\rho_Rv_R-\rho_{C,R}v_C}{\rho_R-\rho_{C,R}}.
$
These four speeds and their associated trajectories partition space into five intervals. The solution is constant inside each interval, except inside the rarefaction wave, where it takes the form:
\begin{equation}
	\begin{array}{ll}
		v_{\text{rare}}(\xi) = \brk*{(\gamma-1)v_L + 2(c_l+\xi)}\prt*{\gamma+1}^{-1},
		 &
		c_{\text{rare}}(\xi) = v_{\text{rare}}(\xi) - \xi,
		\\
		\rho_{\text{rare}}(\xi) = \brk*{\rho_L^\gamma c_{\text{rare}}^2(\xi) \prt*{\gamma \curlyP_L}^{-1}}^{\frac{1}{\gamma-1}},
		 &
		\curlyP_{\text{rare}}(\xi) = \rho_{\text{rare}}^\gamma(\xi) \curlyP_L \rho_L^{-\gamma},
	\end{array}
\end{equation}
where $\xi = x/t$. In summary, at time $t$ and position $x$, the solution takes the form
\begin{equation}
	\left\{
	\begin{array}{llll}
		\displaystyle
		\rho = \rho_L,                        &
		v = v_L,                              &
		\curlyP = \curlyP_L,                  &
		\text{if } \xi \leq v_L-c_L;
		\\\displaystyle
		\rho = \rho_{\text{rare}}(\xi),       &
		v = v_{\text{rare}}(\xi),             &
		\curlyP = \curlyP_{\text{rare}}(\xi), &
		\text{if } v_L-c_L \leq \xi \leq v_C-c_C;
		\\\displaystyle
		\rho = \rho_{C,L},                    &
		v = v_C,                              &
		\curlyP = \curlyP_C,                  &
		\text{if } v_C-c_C \leq \xi \leq v_C;
		\\\displaystyle
		\rho = \rho_{C,R},                    &
		v = v_C,                              &
		\curlyP = \curlyP_C,                  &
		\text{if } v_C \leq \xi \leq s;
		\\\displaystyle
		\rho = \rho_R,                        &
		v = v_R,                              &
		\curlyP = \curlyP_R,                  &
		\text{if } s \leq \xi.
	\end{array}
	\right.
\end{equation}

A two fluid Sod Tube problem can be constructed by considering two initially separate fluids, in contact at the point $x=0$; the left fluid with state ($\rho_L$, $v_L$, $\curlyP_L$), and the right fluid with state ($\rho_R$, $v_R$, $\curlyP_R$). This datum generates a solution profile which is identical to the one described above, and the two fluids simply remain unmixed, separated by the contact discontinuity \cite{LeVeque2002}. This is true even if the fluids have different adiabatic exponents, in which case the solution is calculated by using $\gamma_L$ and $\gamma_R$ respectively left and right of the contact discontinuity.

\FloatBarrier

{
	\small
	\bibliographystyle{abbrv}
	\bibliography{./BailoRey_ProjectiveKineticMixtures.bib}
	\addcontentsline{toc}{section}{References}
}

\end{document}